\documentclass[preprint,12pt]{elsarticle1}

\usepackage{graphicx}
\usepackage{amssymb}
\usepackage{amsthm}
\usepackage{multirow}
\usepackage[english]{babel}
\usepackage{color}

\usepackage{nfontdef}

\newcommand{\ignore}[1]{}

\newcommand{\vek}[1]{\mathchoice{\displaystyle\boldsymbol#1}
{\textstyle\boldsymbol#1}{\scriptstyle\boldsymbol#1}
{\scriptscriptstyle\boldsymbol#1}}
\newcommand{\mat}[1]{\mathchoice{\displaystyle\mathbf#1}
{\textstyle\mathbf#1}{\scriptstyle\mathbf#1}
{\scriptscriptstyle\mathbf#1}}

\newcommand{\di}{\mathrm{d}}

\journal{Applied Mathematics and Computation}

\begin{document}

\begin{frontmatter}

  \title{Uncertainty Updating in the Description of Coupled Heat and
    Moisture Transport in Heterogeneous Materials}

  \author[ctu]{Anna Ku\v{c}erov\'a\corref{auth}}
  \ead{anicka@cml.fsv.cvut.cz}
  \author[ctu,cideas]{Jan S\'{y}kora} \ead{jan.sykora.1@fsv.cvut.cz}
  \cortext[auth]{Corresponding author.
    Tel.:~+420-2-2435-5326;fax~+420-2-2431-0775}
  \address[ctu]{Department of Mechanics, Faculty of Civil Engineering,
    Czech Technical University in Prague, Th\' akurova 7, 166 29
    Prague 6, Czech Republic}
  \address[cideas]{Centre for Integrated Design of Advances
    Structures, Th\' akurova 7, 166 29 Prague 6, Czech Republic}

\begin{abstract}
  To assess the durability of structures, heat and moisture transport
  need to be analyzed. To provide a reliable
  estimation of heat and moisture distribution in a certain structure, one needs to include all
  available information about the loading conditions and material parameters. Moreover, the information should be
  accompanied by a corresponding evaluation of its
  credibility.  Here, the Bayesian inference is applied to combine
  different sources of information, so as to provide
  a more accurate estimation of heat and moisture fields
  \cite{Tarantola:2005}. The procedure is demonstrated
  on the probabilistic description of heterogeneous material where
  the uncertainties consist of a particular value of individual material characteristic
  and spatial fluctuations. As for the heat and
  moisture transfer, it is modelled in coupled setting
  \cite{Kunzel:1997}.
\end{abstract}


\end{frontmatter}

\section{Introduction}
\label{sec:intro}
There are many important factors limiting the service life of
buildings. An appropriate reliability analysis needs to take into
account uncertainties in the environmental conditions as well as
in structural properties. Thanks to the growth of powerful
computing resources and technology, recently developed procedures
in the field of stochastic mechanics have become applicable to
realistic engineering systems.

The most common methods quantifying uncertainties are the first- and
second-order reliability methods (FORM/SORM \cite{Ditlevsen:1996})
computing the probability of failure related to limit states.
Nevertheless, modern sophisticated and highly nonlinear models lead to
non-Gaussian higher-order probability density functions of model
parameters and its response. The higher moments become interesting
quantities to be estimated by probabilistic analysis such as
stochastic finite element methods (SFEM), see e.g.
\cite{Stefanou:2009:CMAME} for a recent review.
SFEM is a powerful tool in computational stochastic mechanics
extending the classical deterministic finite element method (FEM) to the
stochastic framework involving finite elements whose properties are
random.

\begin{figure}[h!]
\centering
\includegraphics[keepaspectratio,width=3.5cm]{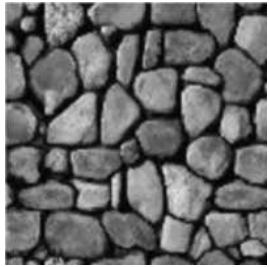}
\caption{Example of quarry masonry}
\label{fig_irrmasonry}
\end{figure}

This paper is focused on the modelling of uncertainties in
material properties and investigates the influence of such
uncertainties on structural response. When dealing with
homogeneous materials, one obtains a simplified scenario for SFEM,
where uncertain material properties are described by random
variables, which are assumed to be spatially constant. In the
field of heterogeneous materials modelling, a widespread approach
is multiscale modelling based on homogenization as presented e.g.
in \cite{Sejnoha:MS:2008,Valenta:2009:IJMCE,Vorel:2008:SEM}.
Nevertheless, homogenization theories have rigorous foundations
for materials with well-defined geometry and components described
by simple constitutive laws as in the case of regular masonry or
composite materials with periodic microstructure. However, these
techniques cannot be efficiently applied to materials with random
microstructure such as in case of quarry masonry,
see~Fig.~\ref{fig_irrmasonry}.

Another possibility is casting the description of heterogenous
material within the probabilistic framework, where uncertain
material properties in time and/or space are represented by
stochastic processes and fields. The resulting problem can be then
solved by different strategies. The most famous group involves
spectral stochastic finite element methods (SSFEM)
\cite{Ghanem:2003}. These methods of uncertainty quantification
are focused on propagation of uncertainty in the system properties
through the numerical model of the system in order to estimate
probabilistic quantification of structural response. It assumes
the knowledge of probabilistic formulation of uncertain system
properties. Nevertheless, this uncertainty is usually very high
before the structure is built, but after the construction,
practical measurements can be performed. With these observations,
the probabilistic models can be updated, to give a more accurate
and reliable estimate. To this goal, the appropriate techniques
from the field of inverse analysis should be employed.

The following section presents a brief introduction into the
inverse analysis. Section \ref{sec:bayes} is devoted to Bayesian
inference suitable for probabilistic estimation of model
parameters from noisy and limited data. Section \ref{sec:het} is
focused on the probabilistic description of heterogeneous
materials properties where particular material parameters are not
spatially constant. Section \ref{sec:example} presents the
application of Bayesian inference to model of coupled heat and
moisture transfer in heterogeneous material and the obtained
results are concluded in Section \ref{sec:concl}.

\section{Inverse problems}\label{sec:inverse}

In computational mechanics one tries to model a real {\it system}
$\msf{A}$, where system {\it parameters} $\msf{q}$, a {\it
loading} $\msf{f}$ and a system {\it response} $\msf{u}$ are
related as
\begin{equation}
\msf{A}(\msf{u};\msf{q}) = \msf{f} \, . \label{eq_system}
\end{equation}
The goal here is to obtain the response of the system for given
parameters and loading conditions. In the field of inverse
analysis, the goal is to find the values of system parameters
$\msf{q}$ corresponding to given loading conditions $\msf{f}$ and
experimental observations $\msf{z}$. Therefore, one uses the
numerical model of the system $\msf{A}$ and derives a so-called
{\it observation
  operator} $\msf{Y}$ mapping the response $\msf{u}$ given parameters $\msf{q}$ and loading $\msf{f}$ to observed
quantities $\msf{z}$
\begin{equation}
\msf{Y}(\msf{q}, \msf{f}) = \msf{z} \, . \label{eq_Adet}
\end{equation}

The subject of this work is concerned with the description of heat
and moisture conduction in structural materials and the system
parameters $\msf{q}$ are related here to the material behaviour.
Material parameters are usually determined in the context of a
chosen experimental setup, where the loading conditions are fixed,
hence, the loading $\msf{f}$ is assumed to be constant in the
following text.

When simulating some real experiment, the model response is
usually not equal to measured data because of experimental errors
or imperfection of numerical model itself. It is often difficult
to distinguish these sources of errors and they are described
together in error vector $\msf{\epsilon}$, so Eq. (\ref{eq_Adet})
becomes
\begin{equation}
\msf{Y}(\msf{q}) + \msf{\epsilon} = \msf{z} \, . \label{eq_Adet01}
\end{equation}

The most common way of estimating material properties is based on
fitting the response of numerical model to the results of real
experiments, see e.g. \cite{Kucerova:2007:PHD} for a recent review
of parameters identification strategies. This usually leads to an
optimization problem, where the difference between the model
response and measured data is minimized by an appropriate
optimization algorithm. Nevertheless, the formulation of the
suitable error function is not always trivial
\cite{Lubineau:2009:CM,Kucerova:2009:WCSMO,Kuraz:2010:JCAM}. The
resulting function is often multimodal, non-smooth or
non-differentiable and some robust optimization algorithm must be
used (see e.g. \cite{Kucerova:2009:EC} for an illustrative example
of an application of evolutionary algorithm to solve such
problem).

There are only several works on including uncertainties into the
curve fitting-based approach to parameters identification as e.g.
in \cite{Lehky:2005}. In general, these approaches omit related
uncertainties in measurements as well as imperfections of the
numerical model and also the preliminary knowledge about the
material parameters coming from their physical meaning (e.g.
Young's modulus must be positive, Poisson's ratio must lay within
the $-1$ and $0.5$, etc).

\section{Bayesian updating of uncertainties}\label{sec:bayes}

Bayesian inference is the statistical inference in which the
experimental observations are not used as the only source
of~information, but they are used to update the preliminary
probabilistic description of system - the so-called prior
information - to give the posterior distribution
\cite{Kennedy:2001:JRSS}. Recall that in realistic applications,
observations are noisy, uncertain and limited in number relative
to the dimension or complexity of the model space. Also, the model
of a system may have limitations on its predictive value because
of its imprecision, filtering or smoothing effects. Taking into
account all pertinent uncertainties, the process of material
properties estimation cannot lead to a single 'optimal' parameter
set, but one has to find a probability distribution of parameters
that represents the~knowledge about parameter values. The Bayesian
setting for the inverse problems offers a rigorous foundation for
inference from noisy data and uncertain forward models, a natural
mechanism for incorporating prior information, and a quantitative
assessment of uncertainty in the inferred results summarizing all
available information about the unknown quantity
\cite{Tarantola:2005}. In addition, unlike other techniques that
aim to regularize the ill-posed inverse problem to achieve a point
estimate, the Bayesian method treats the inverse problem as a
well-posed problem in an expanded stochastic space.

The Bayesian approaches to inverse problems have received much
recent interest, since increasing performance of~modern computers
and clusters enables exhaustive Monte Carlo computations.  Among
recent applications one can cite applications in environmental
modelling \cite{Yee:2008:JWEIA}, hydrology \cite{Fu:2009:JHydro}
or heat transfer \cite{Parthasarathy:2008:IJHMT}. We review this
approach briefly below; for more extensive introductions, see
\cite{Tarantola:2005}.

The main principle of Bayesian inference is casting the inverse
problem in the probabilistic setting, where material parameters
$\msf{q}$ as well as observations $\msf{z}$ and also the response
of forward operator $\msf{Y}(\cdot)$ are considered as~random
variables or random fields. Therefore, we introduce the following
notation. We consider a set $\mit{\Omega}$ of~random elementary
events $\omega$ together with $\sigma$-algebra $\E{S}$ to which a
real number in the interval $[0,1]$ may be assigned,
the~probability of occurrence - mathematically a measure $\D{P}$.

In the Bayesian setting, we assume three sources of information
and uncertainties, which should be taken into account. The first
one is our prior knowledge about the model/material parameters
$\msf{q}(\omega)$, which is represented by defining the prior
density function $p_{\msf{q}}(\msf{q})$. Prior models may embody
simple constraints on $\msf{q}$, such as a range of feasible
values, or may reflect more detailed knowledge about the
parameters, such as correlations or smoothness.

Other source of information comes from measurements, which are
violated by uncertain experimental errors
$\msf{\epsilon}(\bar{\omega})$. Last uncertainty arises from
imperfection of the numerical model included in the observation
operator $\msf{Y}(\cdot)$, when for example our description of the
real system $\msf{A}$ does not include all important phenomena and
therefore the forward operator response $\msf{Y}(\cdot,
\bar{\bar{\omega}})$ can be assumed as uncertain. The
probabilistic formulation of Eq. (\ref{eq_Adet01}) now becomes
\begin{equation}
\msf{z} = \msf{Y}(\msf{q}(\omega),\bar{\bar{\omega}}) +
\msf{\epsilon}(\bar{\omega}) \, .
\end{equation}

If modelling uncertainties $\bar{\bar{\omega}}$ cannot be
neglected, they can be described by conditional probability
density $p(\msf{z}|\msf{q})$ for predicted data $\msf{z}$ and
given model parameters $\msf{q}$. If these uncertainties can be
neglected, only model parameters $\msf{q}(\omega)$ and
observations $\msf{z}(\bar{\omega})$ remain uncertain. In
practise, it is sometimes difficult to distinguish the
imperfection of the system description $\msf{A}$ from measurement
error $\msf{\epsilon}$. Hence modelling uncertainties
$\bar{\bar{\omega}}$ can be hidden in measuring error
$\msf{\epsilon}(\bar{\omega})$. Finally, for noisy measurements we
define the last probability density $p_{\msf{z}}(\msf{z})$.

To update our prior knowledge about model parameters we must
include measurements with our theoretical knowledge. Bayesian
update is based on the idea of Bayes' rule defined for
probabilities. Definition of Bayes' rule for continuous
distribution is, however, more problematic and hence \cite[Chapter
1.5]{Tarantola:2005} derived the posterior state of information
$\pi(\msf{q},\msf{z})$ as a conjunction of all information at hand
\begin{equation}
\pi(\msf{q},\msf{z}) = \kappa p_{\msf{q}}(\msf{q})
p_{\msf{z}}(\msf{z}) p(\msf{z}|\msf{q}) \, , \label{eq_posterior}
\end{equation}
where $\kappa$ is a normalization constant.

The posterior state of information defined in the space of model
parameters $\msf{q}$ is given by the marginal probability density
\begin{equation}
\pi_{\msf{q}}(\msf{q}) = \D{E}_{\bar{\omega}} \left[
\pi(\msf{q},\msf{z}) \right] =
\kappa p_{\msf{q}}(\msf{q}) \int_{\mit{\bar{\Omega}}}
p(\msf{z}|\msf{q}) p_{\msf{z}}(\msf{z}) \, \D{P}(\di \bar{\omega})
= \kappa p_{\msf{q}}(\msf{q}) \, L(\msf{q}), \label{eq_margin-m}
\end{equation}
where $\mit{\bar{\Omega}}$ is a set of~random elementary events
$\bar{\omega}$ and measured data $\msf{z}$ enters through the {\it
likelihood function} $L(\msf{q})$, which gives a measure of how
good a~forward operator $\msf{Y}(\msf{q})$ is in explaining the
data $\msf{z}$. Here, $\D{E}_{\bar{\omega}}$ is the expectation
operator averaging over $\mit{\bar{\Omega}}$.

To keep the presentation of different numerical aspects of
particular methods clear and transparent, we focus here on a quite
common and simple case, where modelling-uncertainties are
neglected and measurement errors are assumed to be Gaussian. Then
the likelihood function takes the form
\begin{equation}
L({\msf{q}}) = \kappa\exp \left( - \frac{1}{2} \left(
\msf{Y}({\msf{q}}) - \msf{z} \right)^{\mrm{T}} C^{-1}_{\mrm{obs}}
\left( \msf{Y}({\msf{q}}) - \msf{z} \right) \right) \, ,
\label{eq_likeli}
\end{equation}
where $C_{\mrm{obs}}$ is a covariance among measurements
$\msf{z}$.

The primary computational challenge is extracting information from
the posterior density $\pi_{\msf{q}}(\msf{q})$
\cite{Tarantola:2005}. Most estimates take the form of integrals
over the posterior, which may be computed with asymptotic methods,
deterministic methods, or sampling. The deterministic quadrature
or cubature may be attractive alternatives to Monte Carlo
simulation at low to moderate dimensions, but Markov chain Monte
Carlo (MCMC) \cite{Tierney:1994:AS} remains the most general and
flexible method for complex and high-dimensional distributions.

\section{Uncertainty in properties of heterogeneous materials}\label{sec:het}

In modelling of heterogeneous material, some material parameters
are not constants, but can be described as random fields. It means
that the uncertainty in the particular material parameter $q$ is
modelled by defining $q(\vek{x})$ for each $\vek{x} \in \C{G}$ as
a random variable $q(\vek{x}): \mit{\Omega} \rightarrow \D{R}$ on
a suitable probability space $(\mit{\Omega} ,\mathscr{S},\D{P})$
in some bounded admissible region $\C{G} \subset \D{R}^d$. As a
consequence, $q : \C{G} \times \mit{\Omega} \rightarrow \D{R}$ is
a random field and one may identify $\mit{\Omega}$ with the set of
all possible values of $q$ or with the space of all real-valued
functions on $\C{G}$. Alternatively, $q(\vek{x},\omega)$ can be
seen as a collection of real-valued random variables indexed by
$\vek{x} \in \C{G}$.

Assuming the random field $q(\vek{x}, \omega)$ to be Gaussian, it
is defined by its mean
\begin{equation}
\mu_q(\vek{x}) = \D{E}[q(\vek{x},\omega)] = \int_{\mit{\Omega}}
q(\vek{x},\omega) \, \D{P}(\di \omega)
\end{equation}
and its covariance
\begin{eqnarray}
C_q(\vek{x},\vek{x}^{\prime}) & = & \D{E}[(q(\vek{x},\omega) -
\mu_q(\vek{x}))
(q(\vek{x}^{\prime},\omega) - \mu_q(\vek{x}^{\prime}))] \nonumber \\
& = & \int_{\mit{\Omega}} (q(\vek{x},\omega) - \mu_q(\vek{x}))
(q(\vek{x}^{\prime},\omega) - \mu_q(\vek{x}^{\prime})) \,
\D{P}(\di \omega) \, .
\end{eqnarray}
Some non-Gaussian fields may be synthesized as -- usually
nonlinear -- functions of Gaussian fields
\cite{Matthies:2007:IB,Rosic:2008:JSSCM}. For instance, most of
the material parameters cannot be negative, hence, the lognormal
random field can be more suitable for their description. Then the
lognormal random field $q(\vek{x},\omega)$ for each material
parameter is obtained by nonlinear transformation of a standard
Gaussian random field $q_g(\vek{x},\omega)$ as, given
in~\cite{Rosic:2008:JSSCM},
\begin{equation}
  q(\vek{x},\omega) = \exp ( \mu_g + \sigma_g q_g(\vek{x},\omega)
  ) \, .
\label{eq_logfield}
\end{equation}
The statistical moments $\mu_g$ and $\sigma_g$ can be obtained
from statistical moments $\mu_q$ and $\sigma_q$ given for
lognormally distributed material property according to following
relations:
\begin{equation}
\sigma_g^2  =  \ln \left( 1 + \left( \frac{\sigma_q}{\mu_q}
\right)^2 \right) \, , \quad \mu_g  =  \ln \mu_q - \frac{1}{2}
\sigma_g^2 \, \label{eq_sigmag}.
\end{equation}

In a computational setting, the random field and the numerical
model must be discretized. If the parameter field $q(\vek{x})$ can
be adequately represented on a finite collection of points
$\{{\vek{x}}^n_{i=1}\} \in \D{R}^{2}$, then we can write both the
prior and posterior densities in terms of $\vek{q} =
(q(\vek{x}_1), \dots, q(\vek{x}_n))$, where $q_i = q(\vek{x}_i)$
are random variables usually correlated among each other. The
vector $\vek{q}$, however, will probably be high-dimensional and
it renders MCMC exploration of the posterior more challenging.
Therefore, the Karhunen$-$Lo\`eve expansion (KLE) can be applied
to dimensionality reduction \cite{Marzouk:2009:JCP}.

KLE is an extremely useful tool for the concise representation of
the stochastic processes. Based on the spectral decomposition of
covariance function $C_q(\vek{x},\vek{x}^{\prime})$ and the
orthogonality of eigenfunctions $\psi_{i}$, the random field
$q(\vek{x},\omega)$ can be written as
\begin{equation}
q(\vek{x},\omega) =
\mu_q(\vek{x})+\sum_{i=0}^{\infty}\sqrt{\varsigma_{i}}\xi_{i}(\omega)\psi_{i}(\vek{x}),
\label{eq_kleserie}
\end{equation}
where
$\vek{\xi}(\omega)=(\ldots,\xi_{i}(\omega),\ldots)^{\mathrm{T}}$
is a set of uncorrelated random variables of zero mean and unit
variance. The spatial KLE functions $\psi_i(\vek{x})$ are the
eigenfunctions of the Fredholm integral equation with the
covariance function as the integral kernel:
\begin{equation}
\int_{\C{G}} C_q(\vek{x},\vek{x}^{\prime}) \psi_i(\vek{x})
\mathrm{d} \vek{x}^{\prime} = \varsigma_i \psi_i(\vek{x}) \, ,
\label{eq_fredholm}
\end{equation}
where $\varsigma_i$ are positive eigenvalues ordered in a
descending order.

Since the covariance is symmetric and positive definite, it can be
expanded in the series
\begin{equation}
C_q(\vek{x},\vek{x}^{\prime}) = \sum_{i=1}^{\infty} \varsigma_i
\psi_i(\vek{x}) \psi_i(\vek{x}^{\prime}) \, . \label{eq_covserie}
\end{equation}

However, computing the eigenfunctions analytically is usually not
feasible. Therefore, one discretizes the covariance spatially
according to chosen grid points (usually corresponding to a finite
element mesh). The resulting covariance matrix $\mat{C}_q$ is
again symmetric and positive definite and Eq. (\ref{eq_fredholm})
becomes symmetric matrix eigenvalue problem, see
\cite{Matthies:2005:CMAME}, where the eigenfunctions
$\psi_i(\vek{x})$ are replaced by eigenvectors $\vek{\psi}_i$.
The~eigenvalue problem is usually solved by a Krylov subspace
method with a sparse matrix approximation. For large eigenvalue
problems, \cite{Khoromskij:2008} proposes the efficient low-rank
and data sparse hierarchical matrix techniques.

For practical implementation, the series (\ref{eq_kleserie}) and
(\ref{eq_covserie}) are truncated after $M$ terms, yielding the
approximations
\begin{eqnarray}
\hat{\vek{q}}(\omega) & \approx & \mu_{q} + \sum_{i=1}^{M}
\sqrt{\varsigma_i}\xi_i(\omega) \vek{\psi}_i \, , \label{eq_mkle}\\
\hat{\mat{C}}_q & \approx & \sum_{i=1}^{M} \varsigma_i
\vek{\psi}_i^{\mathrm{T}} \cdot \vek{\psi}_i \, .
\end{eqnarray}
Such spatial semi-discretization is optimal in the sense that the
mean square error resulting from a truncation after the $M$-th
term is minimized.

\section{Uncertainty updating in coupled heat and moisture tranfer}
\label{sec:example}

This section is devoted to application of previously described
techniques to uncertainty updating in coupled heat and moisture
transport in heterogeneous material with uncertain structure as
quarry masonry. In particular, we employ the model proposed by
K\"unzel \cite{Kunzel:1997} described by the energy balance
equation
\begin{equation}
\frac{\mathrm{d}H}{\mathrm{d}\theta}\frac{\mathrm{d}\theta}{\mathrm{d}t}
 =  \vek{\nabla}^{\mathrm{T}}[\lambda\vek{\nabla} \theta]+h_{v}\vek{\nabla}^{\mathrm{T}}
[\delta_{p}\vek{\nabla}\{\varphi p_{\mathrm{sat}}(\theta)\}]
\label{eq:TE03}
\end{equation}
and the conservation of mass equation
\begin{equation}
\frac{\mathrm{d}w}{\mathrm{d}\varphi}\frac{\mathrm{d}
\varphi}{\mathrm{d}t}  =  \vek{\nabla}^{\mathrm{T}}
[D_{\varphi}\vek{\nabla}\varphi]+\vek{\nabla}^{\mathrm{T}}
[\delta_{p}\vek{\nabla}\{\varphi p_{\mathrm{sat}}(\theta)\}] \, ,
\label{eq:TE04}
\end{equation}
where $\theta$ is the temperature, $\varphi$ stands for the
moisture and $H$, $\lambda$, etc are described below. The
transport coefficients defining the material behaviour are
nonlinear functions of structural responses - the temperature and
moisture fields - and material properties. We briefly recall their
particular relations~\cite{Kunzel:1997}:
\begin{itemize}
\item
Thermal conductivity $[\mathrm{Wm}^{-1}\mathrm{K}^{-1}]$:
\begin{equation}
\lambda = \lambda_0 \left( 1+ \frac{b_{\mathrm{tcs}} w_f (b-1)
\varphi}{\rho_s (b - \varphi)} \right) \, . \label{eq_thermcond}
\end{equation}
\item
Evaporation enthalpy of water $[\mathrm{Jkg}^{-1}]$:
\begin{equation}
h_v = 2.5008 \cdot 10^6 \left( \frac{273.15}{\theta}
\right)^{(0.167+3.67 \cdot 10^{-4} \theta)} \, .
\end{equation}
\item Water vapour permeability
$[\mathrm{kgm}^{-1}\mathrm{s}^{-1}\mathrm{Pa}^{-1}]$:
\begin{equation}
\delta_p = \frac{1.9446 \cdot 10^{-12}}{\mu} \cdot \left( \theta +
273.15 \right)^{0.81} \, .
\end{equation}
\item Water vapour saturation pressure $[\mathrm{Pa}]$:
\begin{equation}
p_{\mathrm{sat}} = 611 \exp \left( \frac{17.08 \theta}{234.18 +
\theta} \right) \, .
\end{equation}
\item
Liquid conduction coefficient $[\mathrm{kgm}^{-1}\mathrm{s}^{-1}]$:
\begin{equation}
D_{\varphi} = 3.8 \frac{a^2}{w_f} \cdot
10^{\frac{3w_f(b-1)\varphi}{(b-\varphi)(w_f-1)}} \cdot
\frac{b(b-1)}{(b-\varphi)^2} \, . \label{eq_watvapsat}
\end{equation}
\item
Total enthalpy of building material $[\mathrm{Jm}^{-3}]$:
\begin{equation}
H = \rho_s c_s \theta \, . \label{eq_enthalpy}
\end{equation}
\end{itemize}
More detailed discussion about transport coefficients can be found
in \cite{Kunzel:1997,Cerny:2009:CMEM}. Some transport coefficients
defined by Eqs. (\ref{eq_thermcond}) - (\ref{eq_enthalpy}) depend
on a subset of the material parameters listed in
Tab.~\ref{tab_params}. The approximation factor $b$ appearing in
Eqs. (\ref{eq_thermcond}) and (\ref{eq_watvapsat}) can be
determined from the relation:
\begin{equation}
b = \frac{0.8(w_{80}-w_f)}{w_{80}-0.8w_f} \, ,
\end{equation}
where $w_{80}$ is the equilibrium water content at
$0.8\,\mathrm{[-]}$ relative humidity. Therefore, $b$ is not
considered as a material parameter, while $w_{80}$ is another
material property to be determined. Finally, there are $8$
material parameters listed in Tab. \ref{tab_params} to be
estimated by updating procedure.
\begin{table}[h!]
\centering
\begin{tabular}{lllcc}
\multicolumn{3}{l}{Parameter} & $\mu_q$ & $\sigma_q$ \\
\hline
$w_{f}$ & $\mathrm{[kgm^{-3}]}$ &  free water saturation & 200 & 40 \\
$w_{\mathrm{80}}$ & $\mathrm{[kgm^{-3}]}$ &  water content at $0.8\,\mathrm{[-]}$ relative humidity & 100 & 10 \\
$\lambda_{\mathrm{0}}$ & $\mathrm{[Wm^{-1}K^{-1}]}$ &  thermal conductivity of dry material & 0.3 & 0.1 \\
$b_{\mathrm{tcs}}$ & $\mathrm{[-]}$ &  thermal conductivity supplement & 10 & 2  \\
$\mu$ & $\mathrm{[-]}$ &  water vapour diffusion resistance factor & 12 & 5 \\
$a$ & $\mathrm{[kgm^{-2}s^{-0.5}]}$ &  water absorption
coefficient & 0.6 & 0.2 \\
$c_s$ & $\mathrm{[Jkg^{-1}K^{-1}]}$ &  specific heat capacity & 900 & 100 \\
$\rho_s$ & $\mathrm{[kgm^{-3}]}$ &  bulk density of building material & 1650 & 50 \\
\hline
\end{tabular}
\caption{Mean values and standard deviations of material
parameters} \label{tab_params}
\end{table}

The presented table also contains the prior information about
material parameters in terms of the mean values $\mu_q$ and the
standard deviations $\sigma_q$. Their particular values are chosen
with regard to values corresponding to materials used in masonry
\cite{Pavlik:SFR:2010}. Other prior information is that all these
material parameters cannot be negative and hence, they could be
considered as lognormally distributed. To describe the uncertainty
about parameters of heterogeneous material, we choose the
covariance kernel of a corresponding random field. Since the
material properties change in the space because of changes in
material components, we assume that the spatial fluctuations of
all parameters are equal. It is probably not the best description
of a real material. Nevertheless, we are convinced that the full
spatial correlation among material properties is more realistic
then full spatial independence. For the sake of simplicity, we do
not study here the case of arbitrarily correlated parameters.
Therefore, we assume the same normalized exponential covariance
kernel for all parameters,
\begin{equation}
C(\vek{x}, \vek{x}^{\prime}) = \exp \left( -\left|
\frac{x_1-x^{\prime}_1}{l_{x_{1}}} \right| - \left|
\frac{x_2-x^{\prime}_2}{l_{x_{2}}} \right| \right) \quad
\vek{x}=(x_{1},x_{2}), \,
\vek{x}^{\prime}=(x_{1}^{\prime},x_{2}^{\prime})\, ,
\label{eq_normkern}
\end{equation}
where $l_{x_{1}}$ and $l_{x_{2}}$ are covariance lengths. We
assume also that the expert is certain about correlation lengths
$l_{x_{1}} = 0.1$~$\mathrm{[m]}$ and $l_{x_{2}} =
0.04$~$\mathrm{[m]}$, but he is not sure about a particular
distribution of phases in material. In practise, the correlation
lengths can be determined by the image analysis of a given
material \cite{Lombardo:2009:IJMCE}, common size of bricks in
masonry etc.

Utilizing the covariance kernel (\ref{eq_normkern}), we compute
particular realizations of standard Gaussian random field based on
Karhunen-Lo\`eve expansion
\begin{equation}
\hat{\vek{q}}_g(\omega) \approx \sum_{i=1}^{M}
\sqrt{\varsigma_i}\xi_i(\omega) \vek{\psi}_i \, .
\end{equation}
Here, the eigenvectors $\vek{\psi}_i$ describe the fluctuation of
material property within the studied domain $\C{G}$. Since the
fluctuations reflect here the distribution of phases in material,
the random variables $\xi_i$ are also considered same for all
material parameters. Given the prior mean $\mu_q$ and standard
deviation $\sigma_q$ for each material parameter in Tab.
\ref{tab_params}, corresponding statistical moments $\mu_g$ and
$\sigma_g$ can be derived from Eq. (\ref{eq_sigmag}). These can be
then applied to Eq. (\ref{eq_logfield}) in order to obtain a
lognormal random field for each material parameter. Nevertheless,
such procedure implies that the expert is also certain about mean
value and relative amplitudes (given by the prior standard
deviation) of respective random fields. That will unlikely happen
in practice. To make the example more realistic, we include the
uncertainty in the mean values of particular random fields by
adding one random variable $\xi_{q,0}$ for each material property
and extending the Eq. (\ref{eq_logfield}) into
\begin{equation}
  \hat{\vek{q}}(\omega) = \exp \left( \mu_g + \sigma_g \xi_{q,0} + \sigma_g
\sum_{i=1}^{M} \sqrt{\varsigma_i}\xi_i(\omega) \vek{\psi}_i
  \right) \, .
\label{eq_logfield_ext}
\end{equation}
As a result, the random field corresponding to particular material
property can be now shifted independently to each other. For a
sake of simplicity and to keep the number of random variables in
reasonable bounds, we keep the amplitudes of respective random
fields still closely related. The total number of random variables
$L$ to be updated within Bayesian framework now becomes $L = M+W$,
where $M$ is the number of terms in the truncated KL expansion and
$W$ is the number of material properties.
\begin{figure}[b!]
\centering
\includegraphics[keepaspectratio,width=10cm]{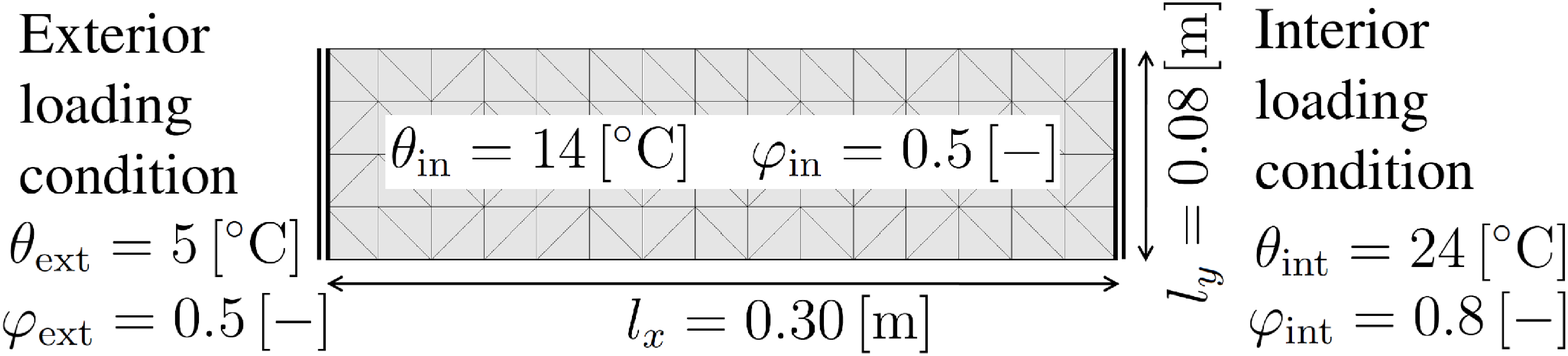}
\caption{Experimental setup} \label{fig_scheme}
\end{figure}

As an example, we consider two-dimensional rectangular domain
discretized by FE mesh into $80$ nodes and $120$ elements. Its
geometry together with the specific loading conditions are shown
in Fig. \ref{fig_scheme}. The initial temperature is
$\theta_{\mathrm{in}} = 14$ [$^{\circ}$C] and moisture
$\varphi_{\mathrm{in}} = 0.5$ [-] in the whole domain. One side of
the domain is submitted to exterior loading conditions
$\theta_{\mathrm{ext}} = 5$ [$^{\circ}$C] and
$\varphi_{\mathrm{ext}} = 0.5$ [-], while the opposite side is
submitted to interior loading conditions $\theta_{\mathrm{int}} =
24$ [$^{\circ}$C] and $\varphi_{\mathrm{int}} = 0.8$ [-].

\begin{figure} [ht!]
\centering
\begin{tabular}{cc}
\includegraphics*[width=65mm,keepaspectratio]{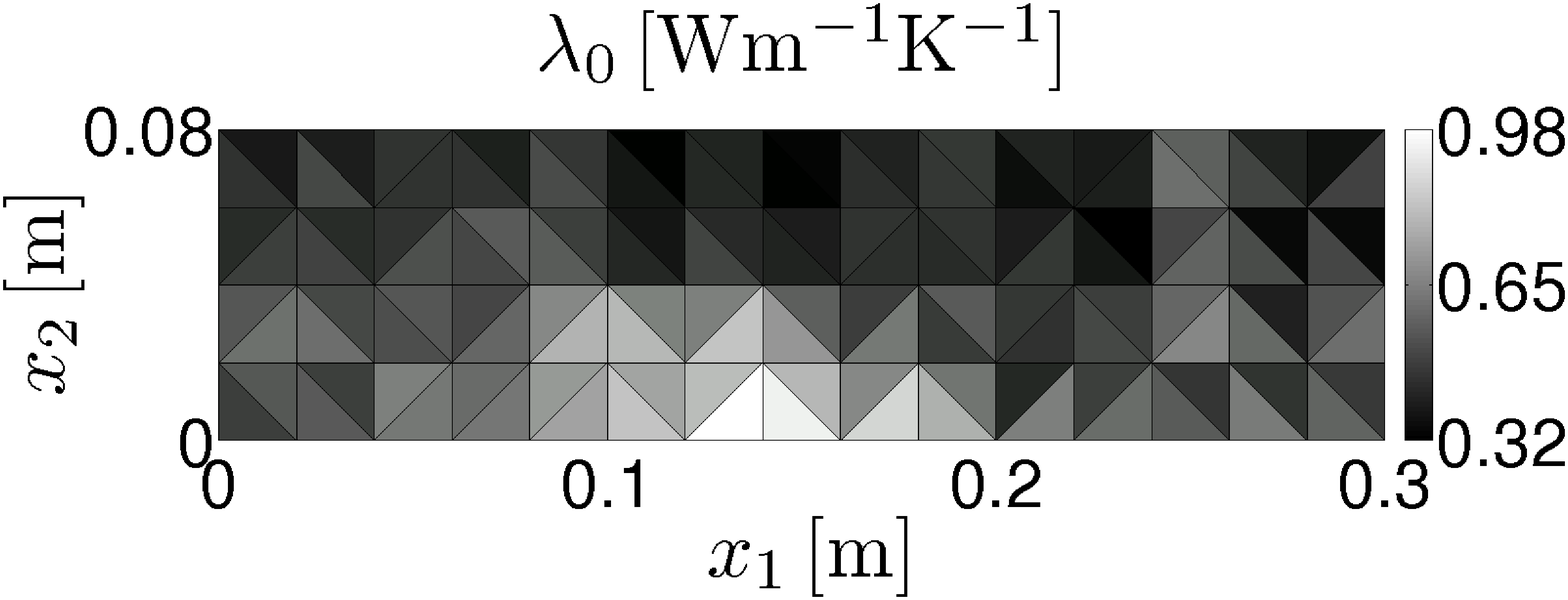}&
\includegraphics*[width=65mm,keepaspectratio]{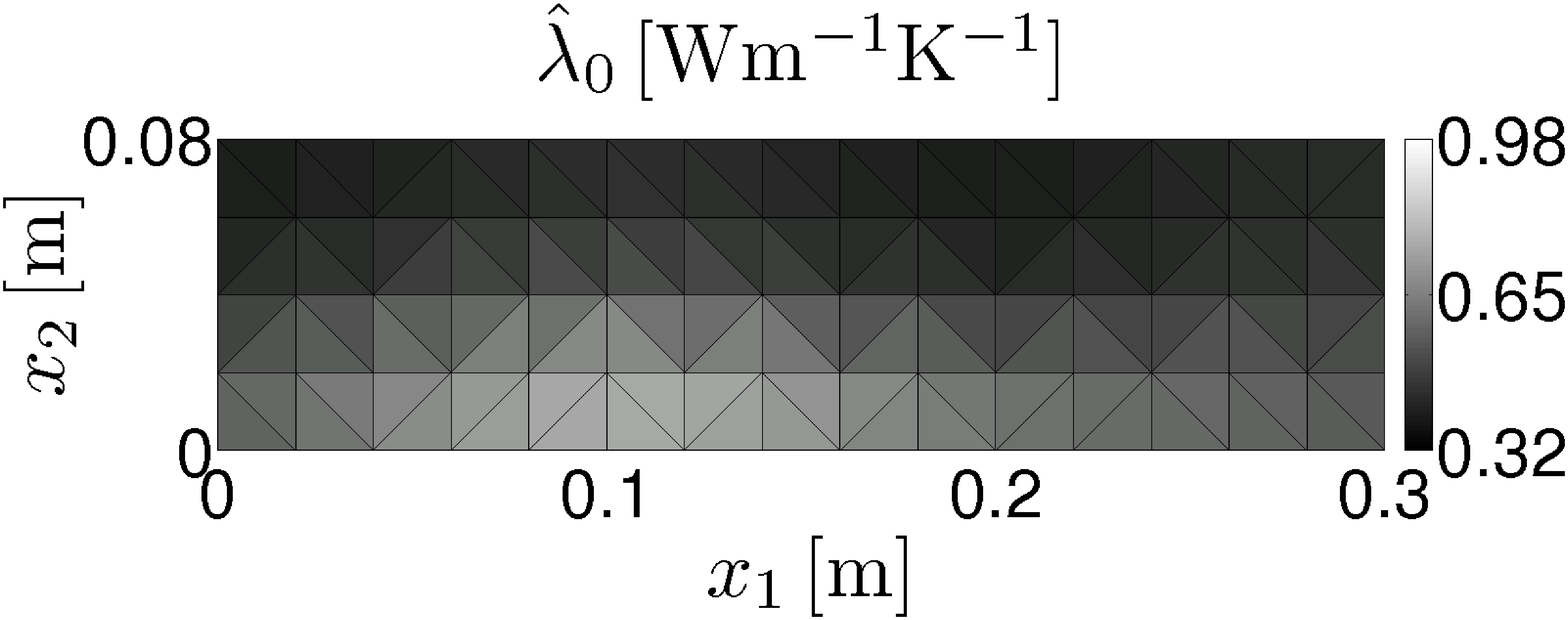}\\
(a)&(b)
\end{tabular}
\caption{Thermal conductivity field computed using (a) all $120$
  eigenmodes and (b) only first $7$ eigenmodes}
\label{fig_lambdafield}
\end{figure}

In order to reduce the number of random variables to be updated
within the Bayesian inference, one needs to choose the number of
eigenmodes $M$ as small as possible, but high enough for a
satisfactory description of parameter fields. The attention should
be paid to the error in description of parameter fields as well as
to the related error in the model response. Fig.
\ref{fig_lambdafield} presents a comparison of an arbitrary
realization of thermal conductivity field $\lambda_0(\vek{x})$
computed using all $120$ eigenmodes and its approximation
$\hat{\lambda}_0(\vek{x})$ computed using only first $7$
eigenmodes.

\begin{figure} [ht!]
\centering
\begin{tabular}{cc}
\includegraphics*[width=65mm,keepaspectratio]{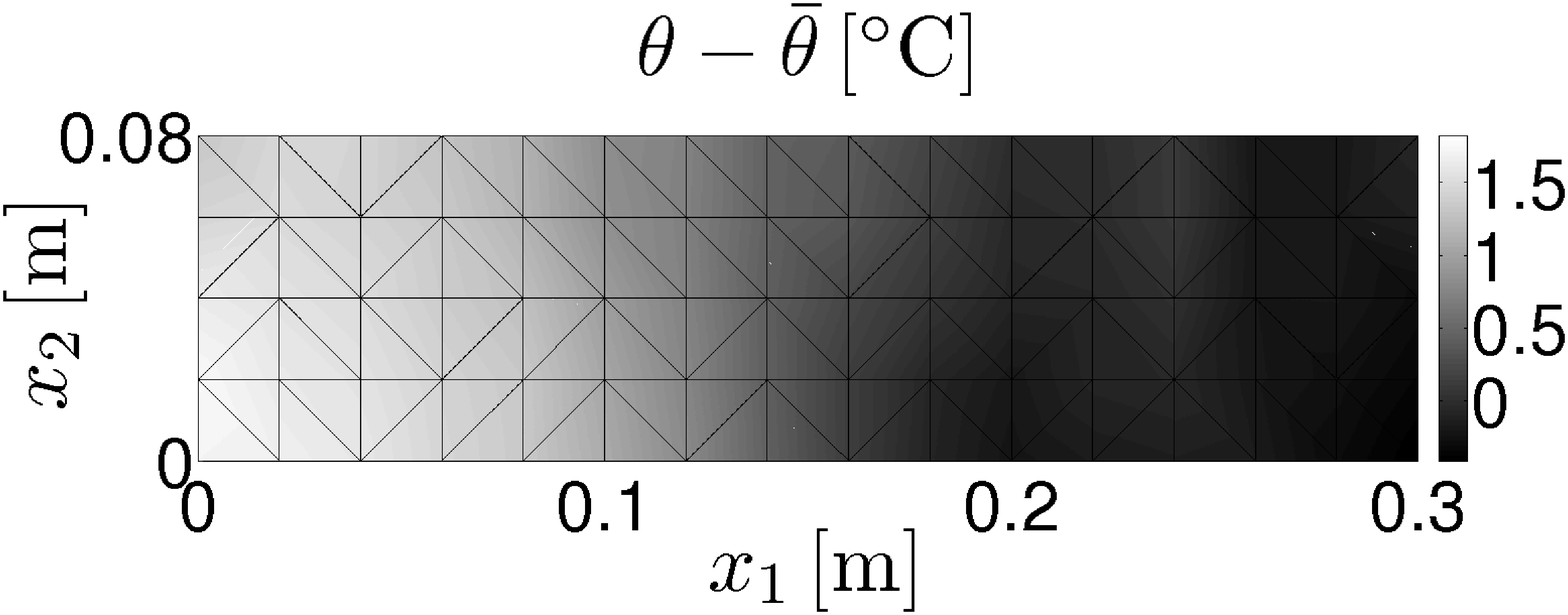}&
\includegraphics*[width=70mm,keepaspectratio]{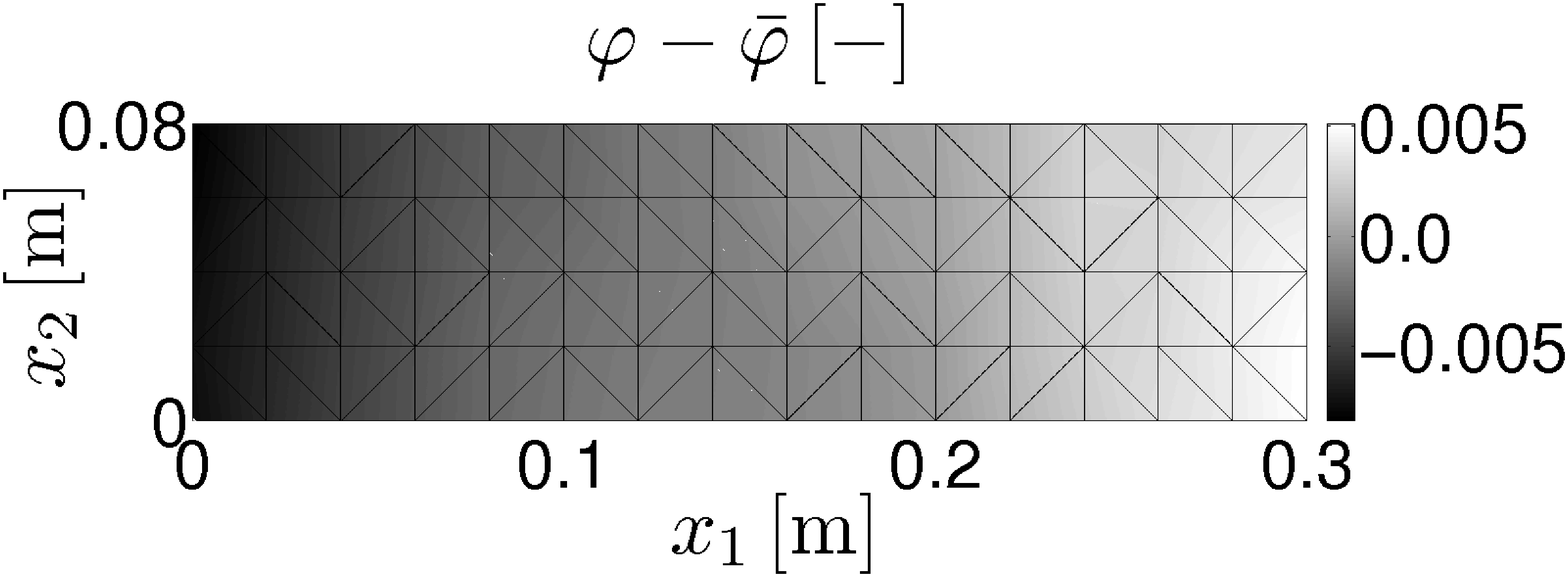}\\
(a)&(b)\\
\includegraphics*[width=65mm,keepaspectratio]{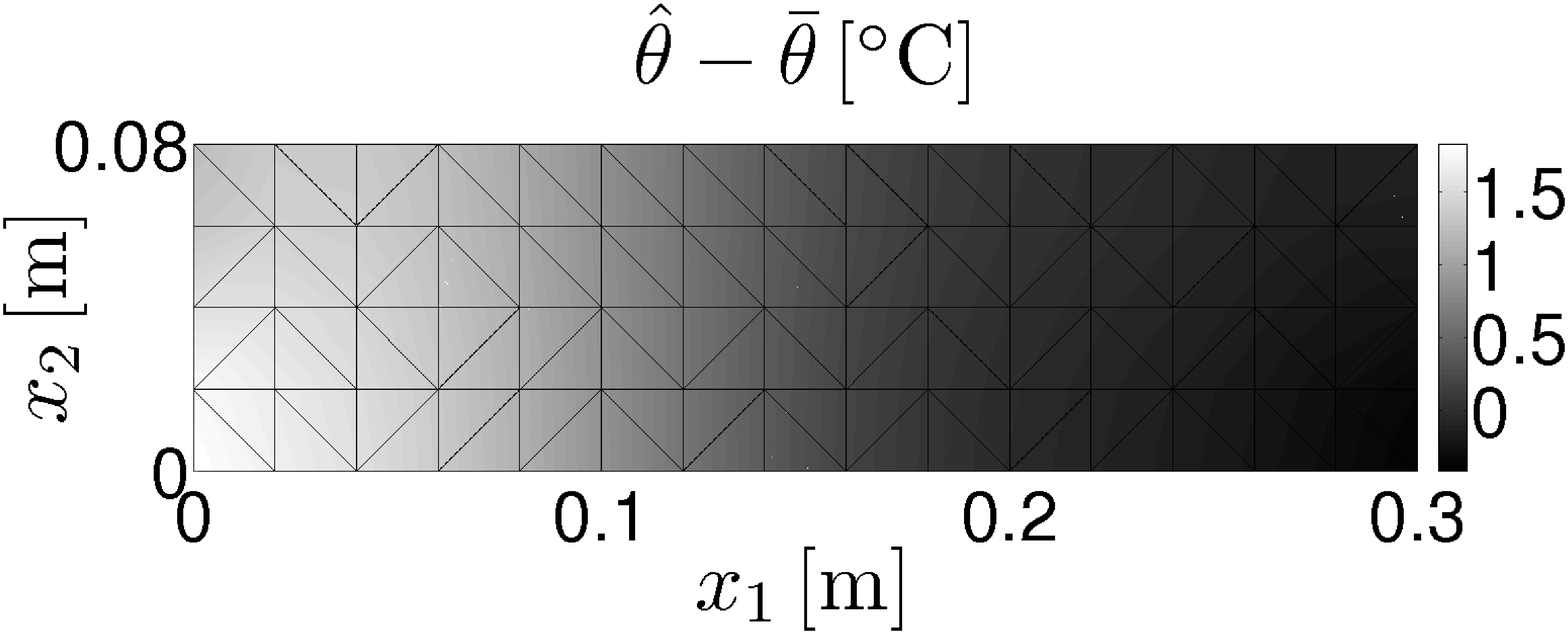}&
\includegraphics*[width=70mm,keepaspectratio]{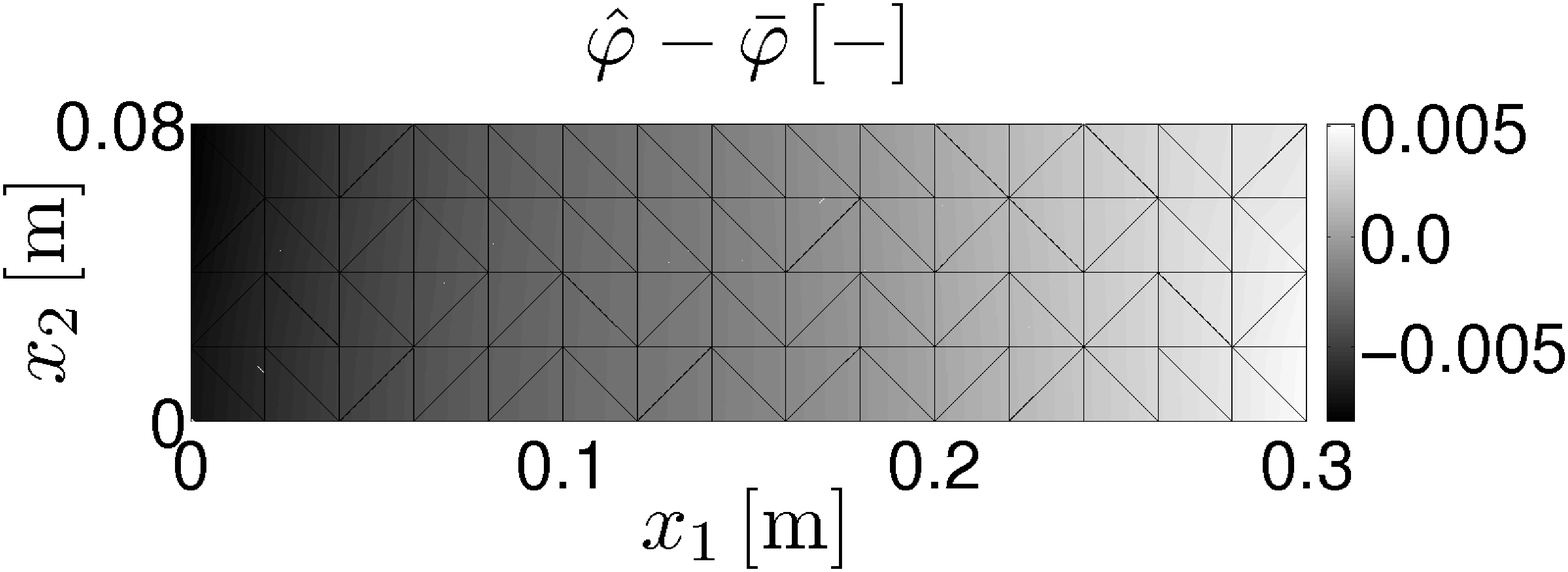}\\
(c)&(d)
\end{tabular}
\caption{Fluctuations of response fields: (a) temperature field and
  (b) moisture field obtained for complete input parameter fields;
  (c) temperature field and (d) moisture field obtained for
  approximated input fields}
\label{fig_respfields}
\end{figure}

\begin{figure} [ht!]
\centering
\begin{tabular}{cc}
\includegraphics*[width=65mm,keepaspectratio]{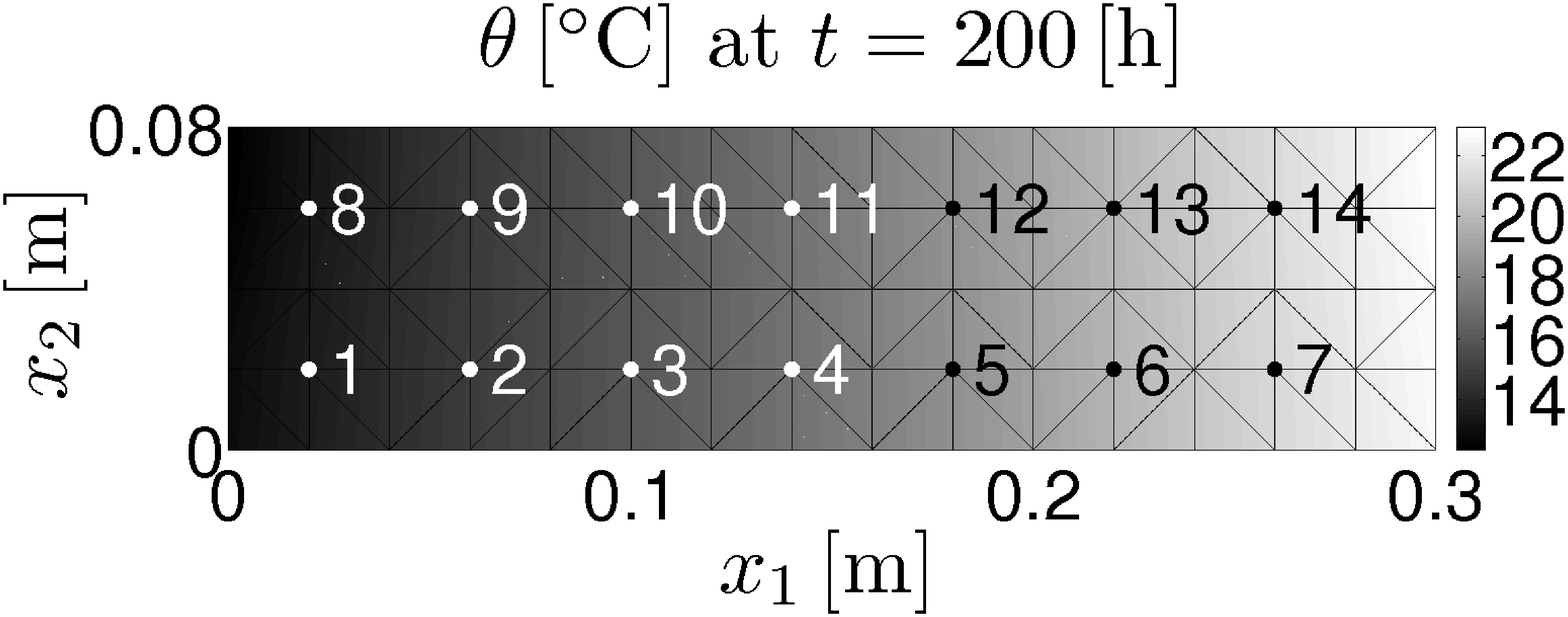}&
\includegraphics*[width=65mm,keepaspectratio]{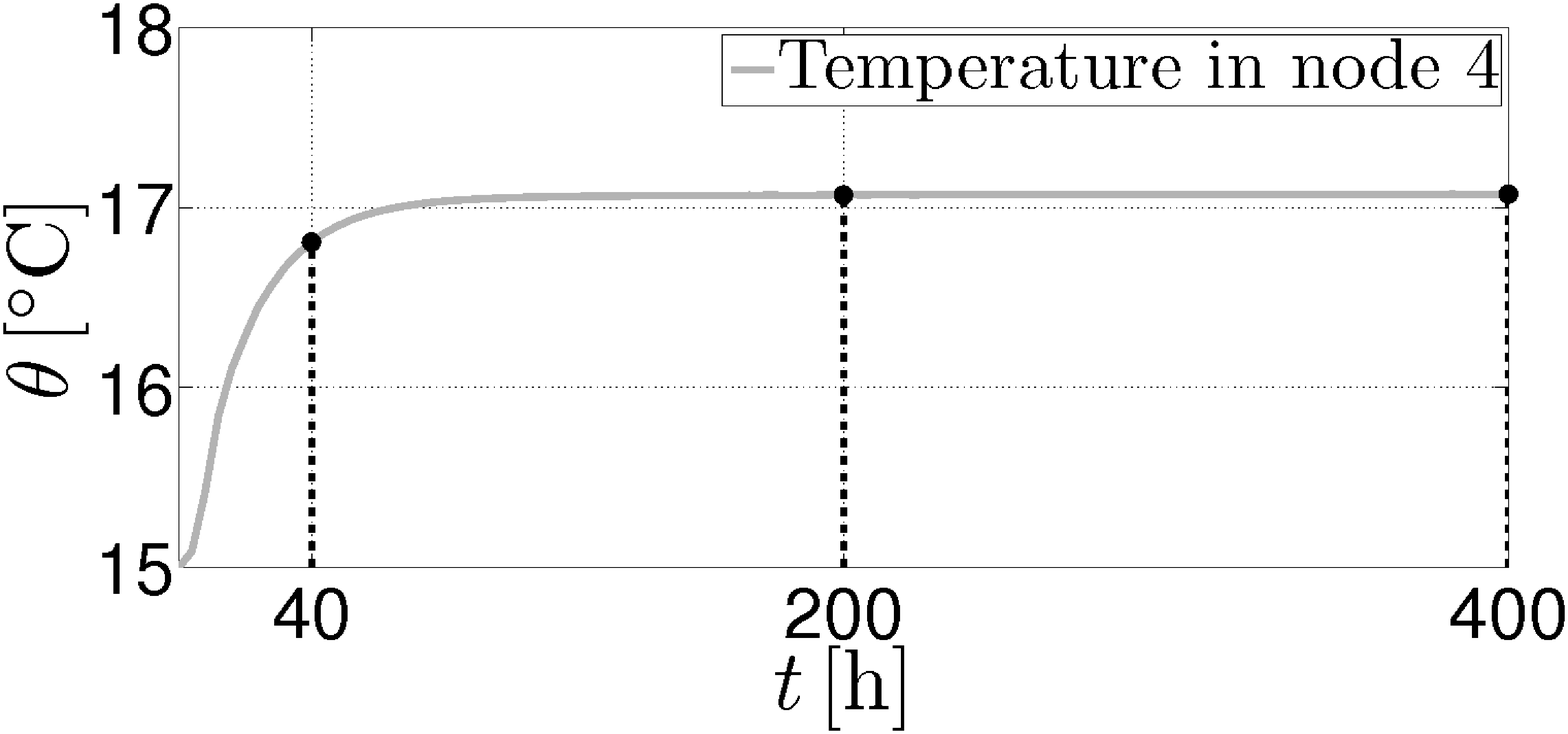}\\
(a)&(b)\\
\includegraphics*[width=65mm,keepaspectratio]{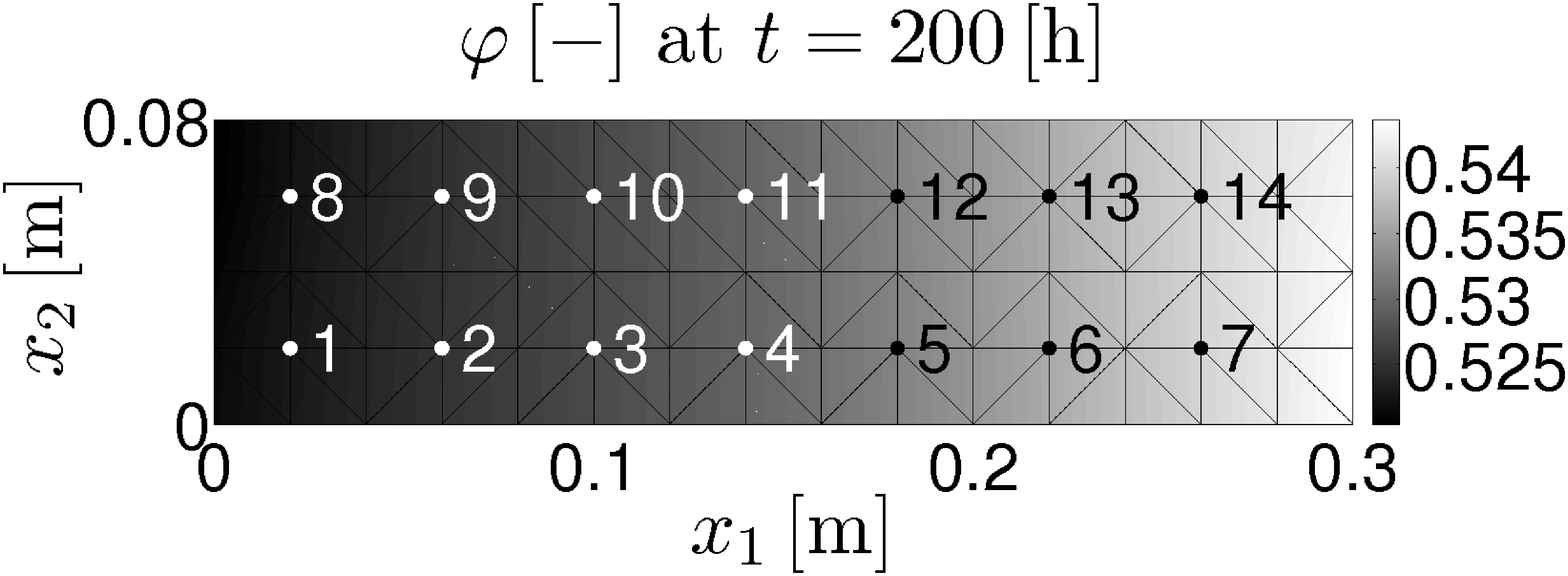}&
\includegraphics*[width=65mm,keepaspectratio]{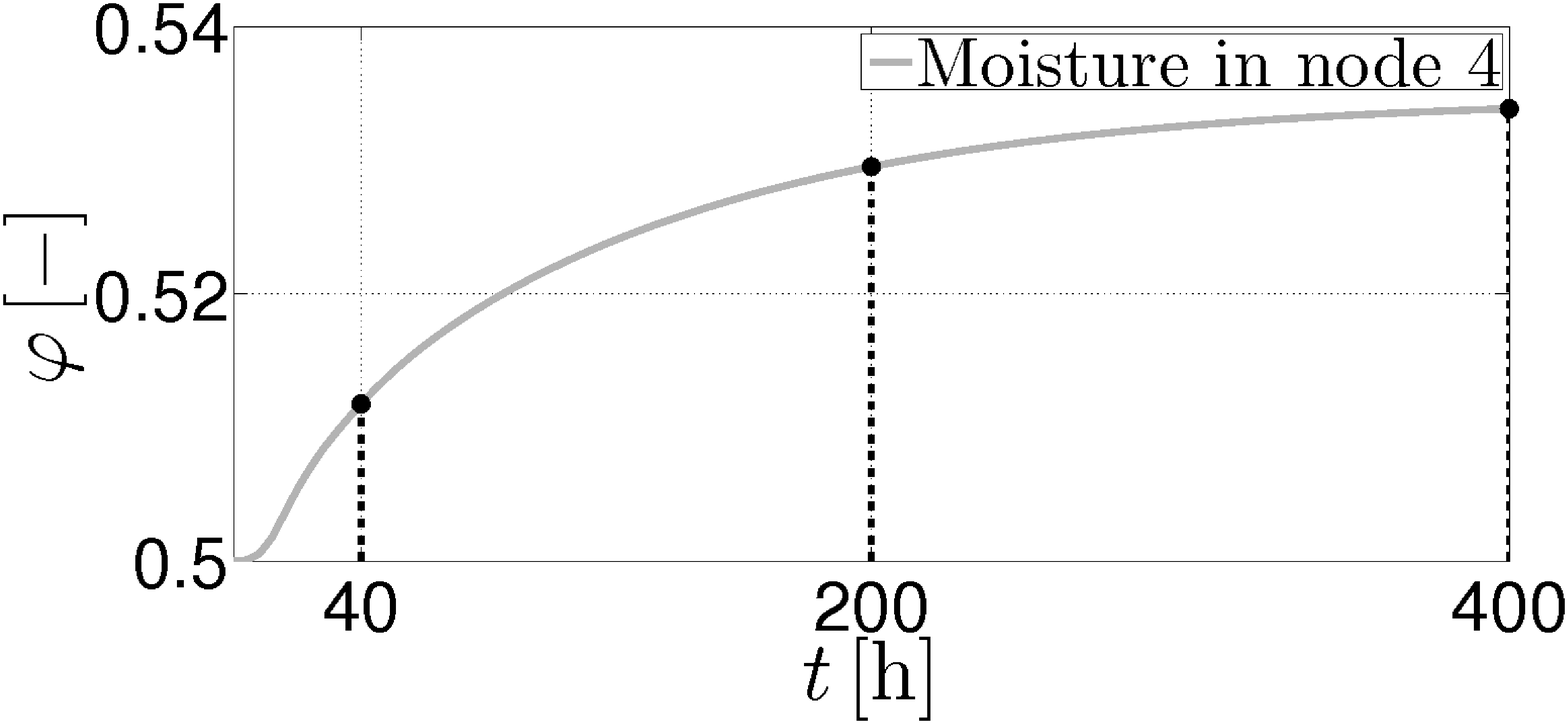}\\
(c)&(d)
\end{tabular}
\caption{Virtual observations: (a) and (c) spatial arrangement of
  probes; (b) and (d) temporal organization of measurements}
\label{fig_observations}
\end{figure}

When comparing the corresponding heat and moisture fields, the
negligence of higher eigenmodes in description of input parameter
fields is reflected by fluctuations of response fields. These are,
however, relatively small comparing to absolute values of
temperature and moisture (see Fig. \ref{fig_observations}). For
better understanding, Fig. \ref{fig_respfields} shows only the
fluctuations of response fields, since we have subtracted the
response fields ($\bar{\theta}(\vek{x})$ and
$\bar{\varphi}(\vek{x})$) corresponding to homogeneous medium.
(The particular response fields correspond to time $t = 200$
$\mathrm{[h]}$.)

One can conclude that the employed K\"unzel's model has a
smoothing effect, because negligence of higher eigenmodes induced
relatively high error in the approximation of input parameter
fields, while its impact to the fluctuations of response fields is
almost vanishing.

\begin{figure} [ht!]
\centering
\begin{tabular}{cc}
\includegraphics*[width=68mm,keepaspectratio]{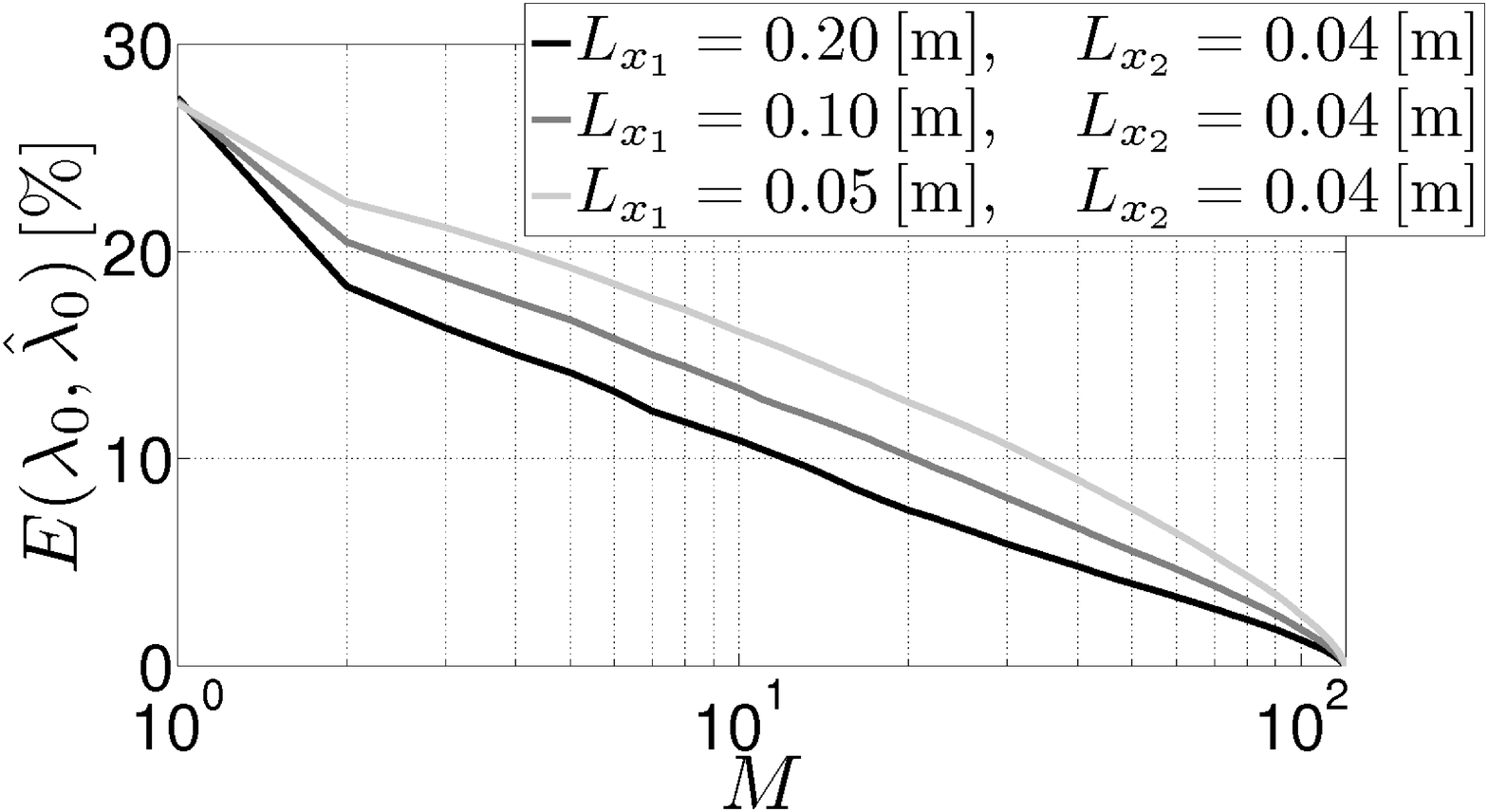}&
\includegraphics*[width=62mm,keepaspectratio]{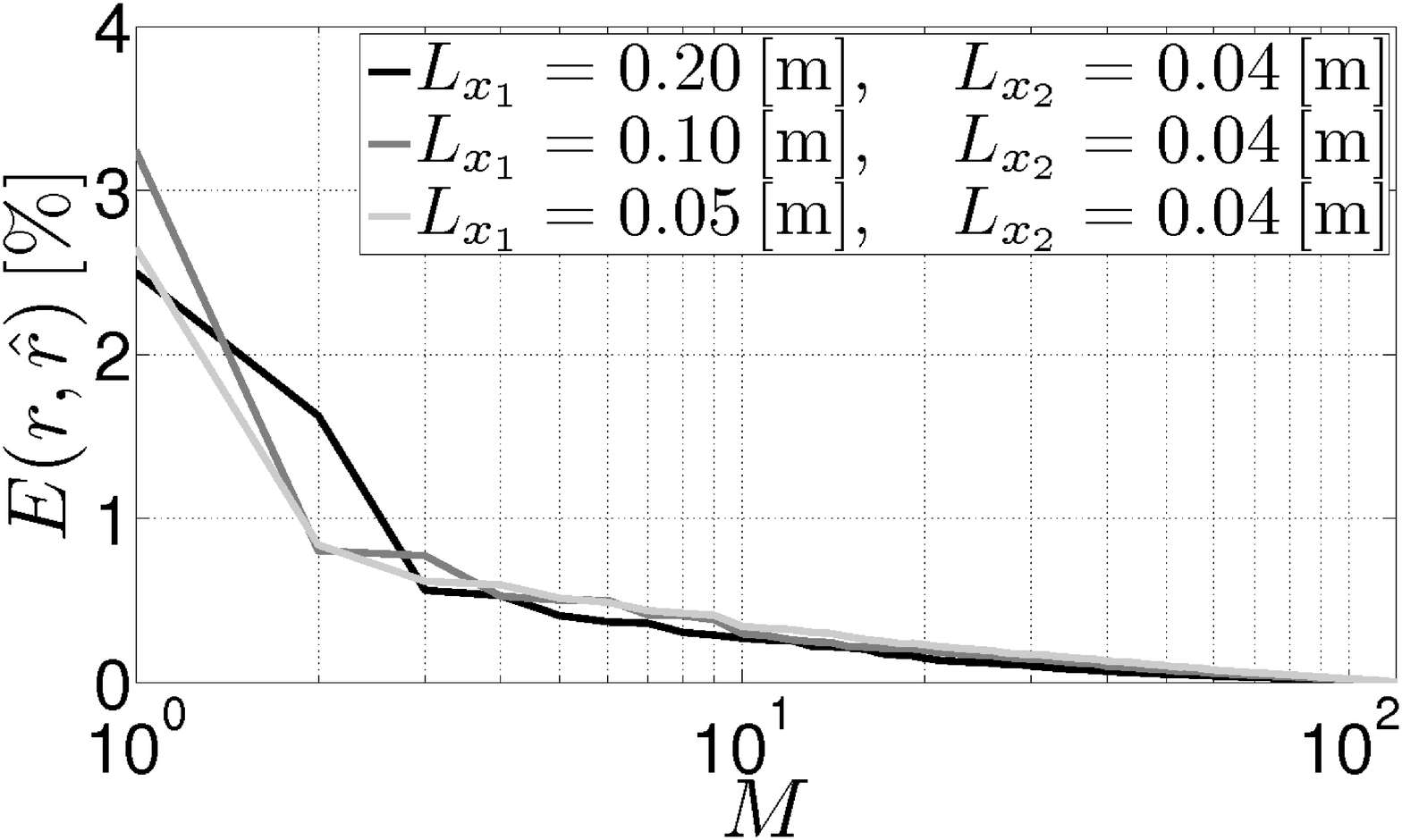}\\
(a)&(b)
\end{tabular}
\caption{Relative mean point-wise error (a) of the input thermal
  conductivity field and (b) of the overall responses induced by KLE
  approximation based on $M$ eigenmodes}
\label{fig_nchierr}
\end{figure}

In order to choose an appropriate number of eigenmodes, a relative
point-wise error of input fields averaged over all $120$ finite
elements and over $100$ independent random realizations can be
computed according to
\begin{equation}
  E(\vek{q},\hat{\vek{q}}) = \frac{1}{100} \sum_{j=1}^{100}
  \frac{1}{120} \sum_{i=1}^{120} \frac{|q_i(\vek{\xi_j}) -
  \hat{q}_i^{(M)}(\vek{\xi_j})|}{q_i(\vek{\xi_j})} \, .
\label{eq_lambdaerr}
\end{equation}
A similar error can be also computed in terms of response fields.
These errors as a function of the number of eigenmodes $M$
involved in the description of input fields computed for three
different choices of correlation lengths are depicted in Fig.
\ref{fig_nchierr}. It can be seen again that the error in
description of input fields is decreasing slowly, while the error
in the response fields descends much faster due to the smoothing
effect of the numerical model. Owing to the results in Fig.
\ref{fig_nchierr}, we decided to perform the Bayesian inference
with an approximation of input fields using $7$ KLE modes, which
provides a satisfactory accuracy.

According to the formulation of lognormal input fields given in
Eq. (\ref{eq_logfield_ext}), we use one random variable for each
eigenmode involved (it is $M = 7$) and one random variable for
each material properties (it is $W = 8$) in order to enable their
relative shift. Hence, we have $L = 15$ random variables to be
updated within the Bayesian inference.

Due to the lack of experimental data, we prepared a virtual
experiment based on simulation including all eigenmodes. The
values of temperature and moisture are measured in 14 points shown
at Figs. \ref{fig_observations} (a) and (c) and at three distinct
times given in Figs. \ref{fig_observations} (b) and (d). Hence,
the observations $\vek{d}$ consist of $84$ values. They were then
perturbed by Gaussian noise with standard deviation for
temperature $\sigma_{\theta} = 0.2$ $\mathrm{[^{\circ}C]}$ and for
moisture $\sigma_{\varphi} = 0.02$ $\mathrm{[-]}$ and in that way
we produced $100$ virtual measurements. Based on them we
calculated the observation covariance matrix
$\mat{C}_{\mathrm{obs}}$ appearing in the likelihood function,
which has the Gaussian form shown in the Eq. (\ref{eq_likeli}).

\begin{figure} [ht!]
\centering
\begin{tabular}{cc}
\includegraphics*[width=65mm,keepaspectratio]{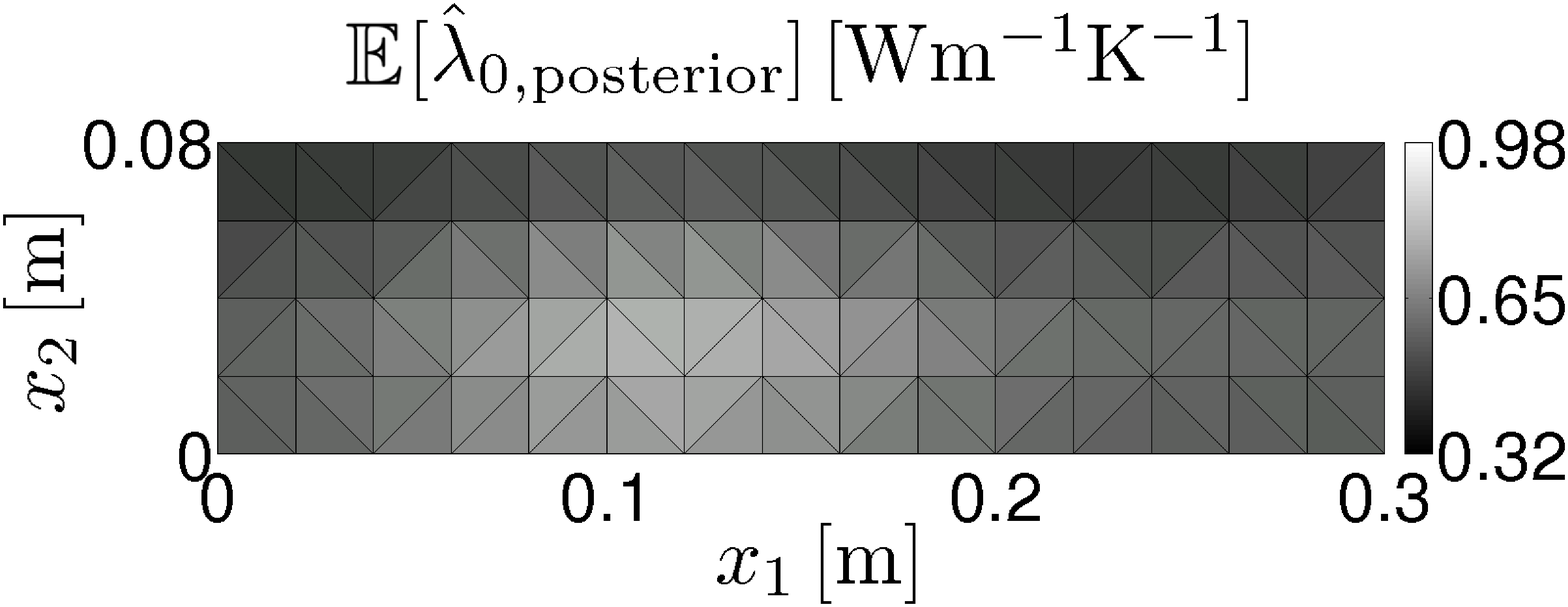}&
\includegraphics*[width=65mm,keepaspectratio]{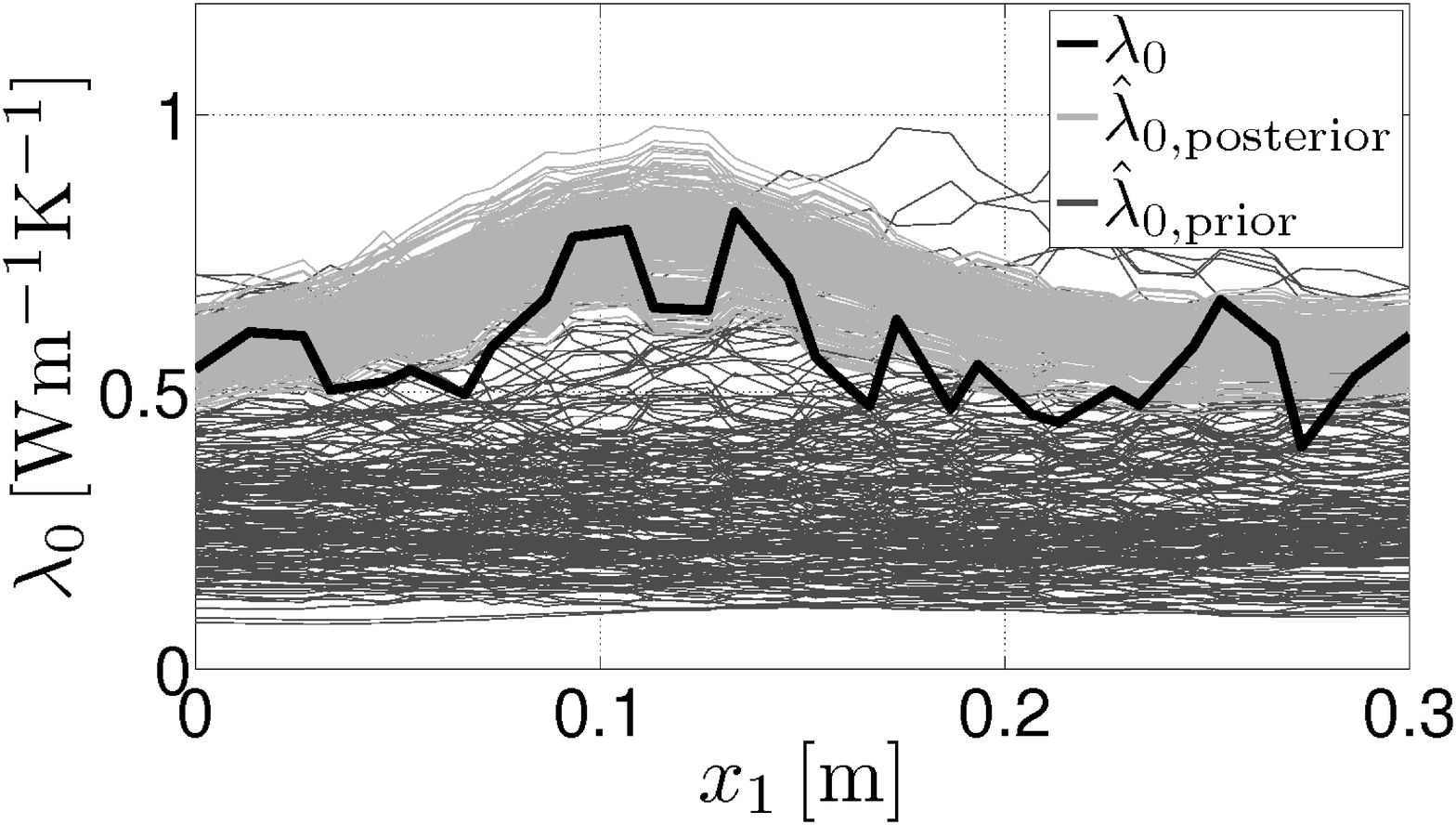}\\
(a)&(b)
\end{tabular}
\caption{A posteriori fields of the thermal conductivity (a) mean over
  $80,000$ samples and (b) cut along the axis $x$ at $y=0.03$ $\mathrm{[m]}$ of a
  subset of $200$ posterior (light grey) and $200$ prior samples (dark
  grey) and the reference field (bold black)}
\label{fig_lambdapostfields}
\end{figure}

The Bayesian update was performed using Metropolis-Hasting
algorithm and $80,000$ samples were generated in order to sample
the posterior density (\ref{eq_posterior}) over variables
$\vek{\xi} = (\xi_1 \dots \xi_{15})$. Then, one realization of
parameter fields was computed for each sample vector $\vek{\xi}$.
The mean computed over all posterior fields of thermal
conductivity $\mathbb{E}[\hat{\lambda}_{0,\mathrm{posterior}}]$
depicted in Fig. \ref{fig_lambdapostfields} (a) can be compared
with the reference field in Fig. \ref{fig_lambdafield} (a)
utilized for preparation of virtual experiment.
Since the mean is not a sufficient descriptor of the a~posteriori
distribution, Fig. \ref{fig_lambdapostfields} (b) shows a cut of a
subset of $200$ posterior samples (light grey lines).  The cut is
driven along the axis $x$ at $y = 0.03$ $\mathrm{[m]}$. These
posterior samples can be compared with the same number of fields
computed for samples drawn from the a priori distribution (dark
grey lines) and the cut through the reference field of $\lambda_0$
(bold black line). One can see that the a posteriori samples much
better encompass the reference field than the a priori ones.

\begin{figure} [ht!]
\centering
\begin{tabular}{cc}
\includegraphics*[width=63mm,keepaspectratio]{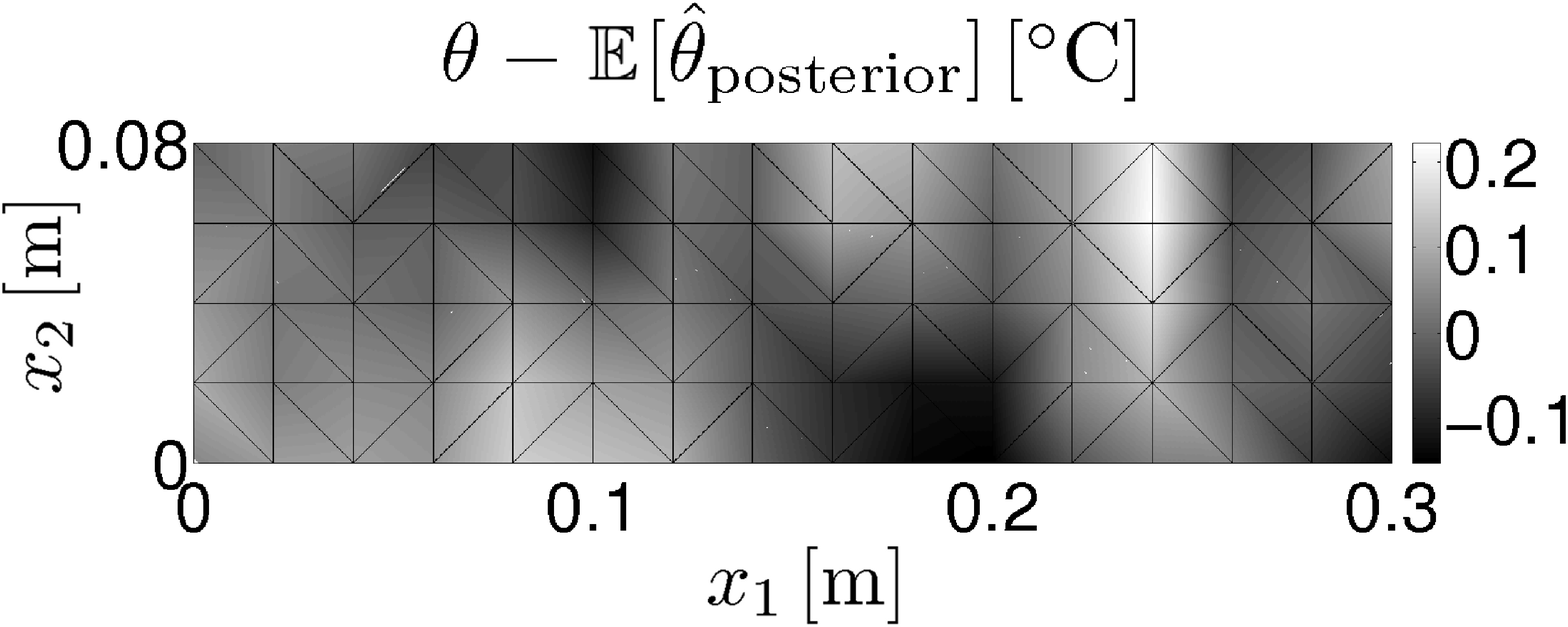}&
\includegraphics*[width=65mm,keepaspectratio]{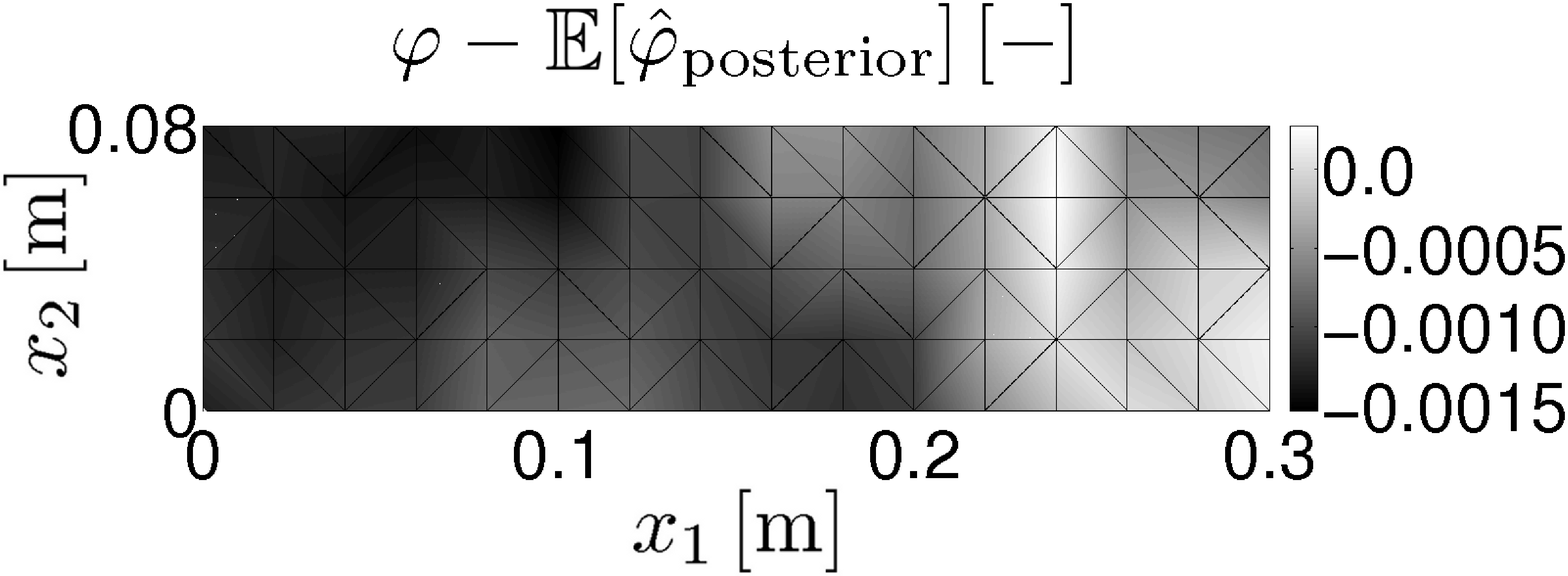}\\
(a)&(b)
\end{tabular}
\caption{Difference between reference field and the mean over the a
  posteriori samples (a) for temperature field and (b) for moisture field}
\label{fig_resppostfields}
\end{figure}

During the updating process, the model responses obtained for the
a~posteriori samples were stored. Fig. \ref{fig_resppostfields}
shows the difference between the reference response fields and the
mean computed over posterior response fields, both at time $t =
200$ $\mathrm{[h]}$. Comparing the differences in Fig.
\ref{fig_resppostfields} with size of fluctuations of
corresponding fields in Figs. \ref{fig_respfields} (a) and (b),
one can conclude that the differences are relatively small.

\begin{figure} [ht!]
\centering
\begin{tabular}{cc}
\includegraphics*[width=80mm,keepaspectratio]{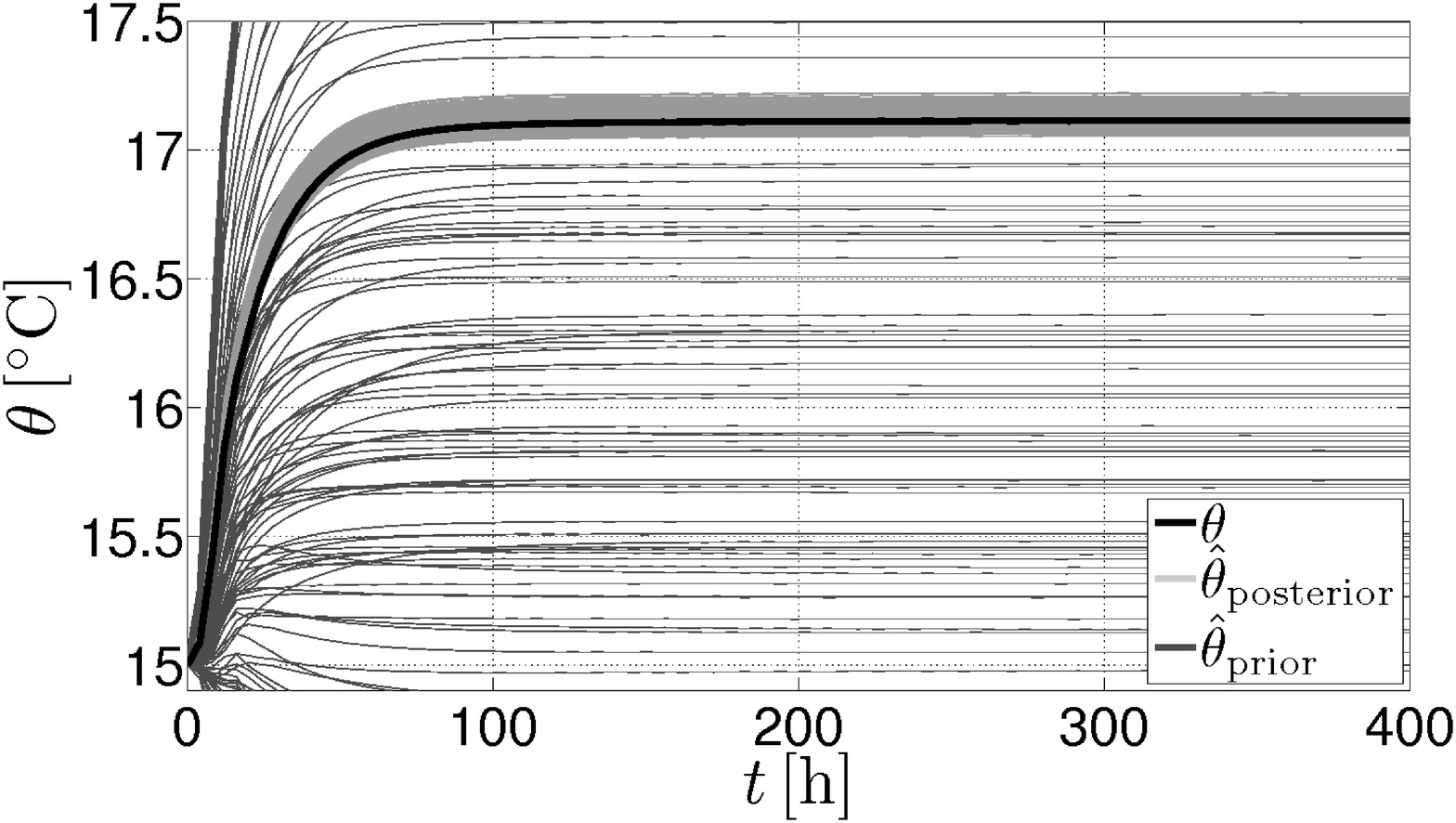}\\
(a)\\
\includegraphics*[width=80mm,keepaspectratio]{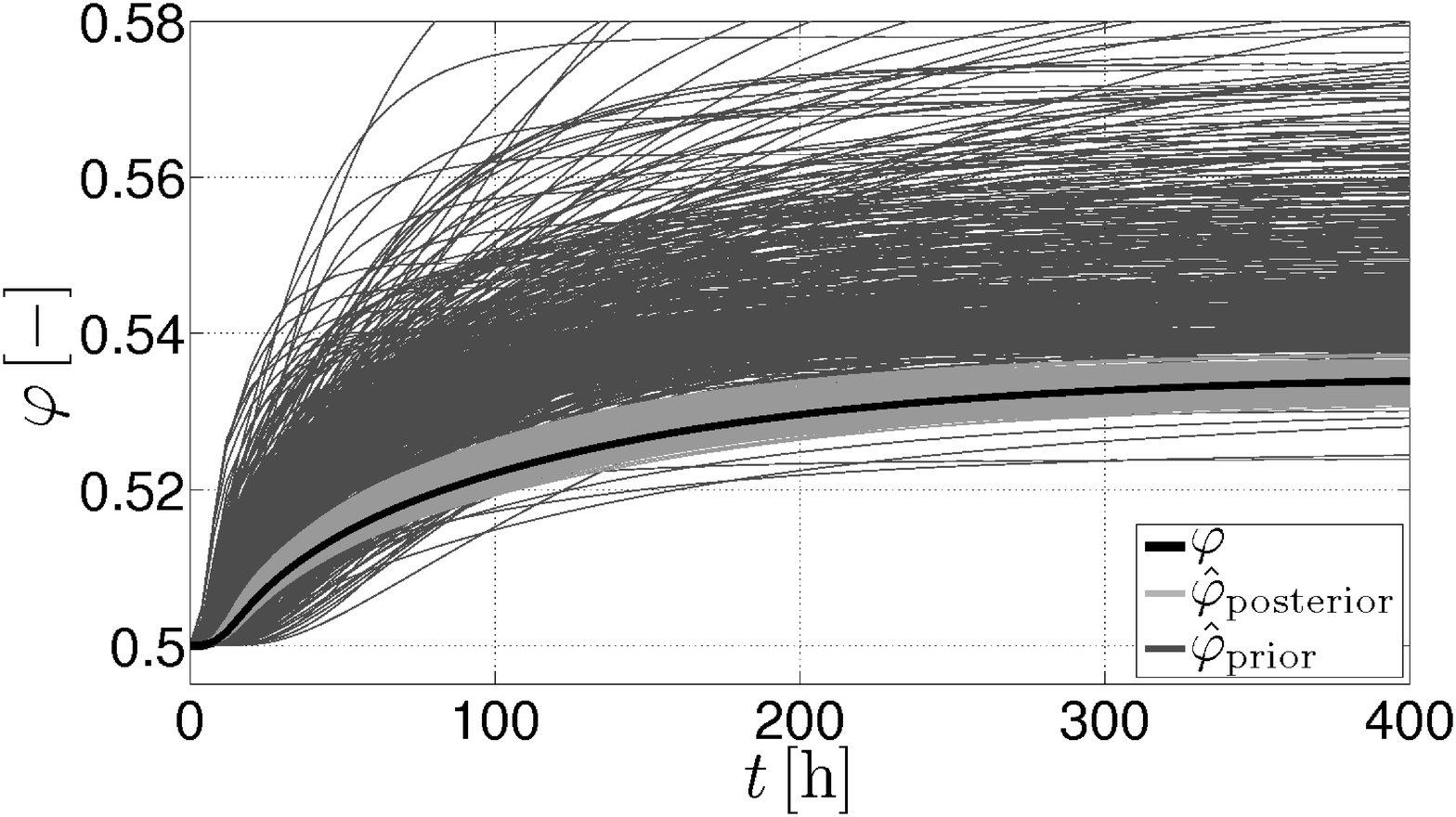}\\
(b)
\end{tabular}
\caption{Temperature (a) and moisture (b) evolution in time at FE node
  $4$: comparison of a priori samples (dark grey), a posteriori
  samples (light grey) and reference response (bold black)}
\label{fig_respposttime}
\end{figure}

Finally, the a posteriori and the a priori distributions of the
responses can be compared with the reference response at Fig.
\ref{fig_respposttime}. The figure shows the evolution of the
temperature (a) and the moisture (b) in time at the FE node no.
$4$. It is evident from the figures that the a posteriori samples
(light grey) are distributed in the very vicinity of the reference
response (bold black), while the dispersion of responses computed
for $200$ samples drawn from the a priori distribution is very
large.

\section{Conclusions}\label{sec:concl}
The presented paper deals with the Bayesian updating of
uncertainty in properties of heterogeneous materials. The process
starts with the a priori highly uncertain information given by an
expert about material properties, which enters as input into the
chosen material model.  The model employed here is the K\"unzel's
model of coupled heat and moisture transfer, which assumes highly
nonlinear relation between $8$ material parameters and the
structural response as given in Eqs. (\ref{eq_thermcond}) --
(\ref{eq_enthalpy}). Beside the prior uncertainty in the value of
material characteristic itself, another uncertainty arises when
describing heterogeneous material with spatial fluctuations of
such property. In order to describe the fluctuations, one random
field is assigned to each material property. For a sake of
simplicity, the fully correlated random fields are assumed, but
one random variable is added into the description of each
parameter so as to allow their mutual shift. To limit the number
of random variables describing the material, the random fields are
approximated by Karhunen-Lo\`eve expansion with 7 eigenmodes.
Finally, to grasp all mentioned uncertainties, $15$ random
variables are considered.

Beside the prior information, the virtual experiment is prepared
to substitute real experimental observations. The
Metropolis-Hastings algorithm is then employed to sample the
posterior distributions of random variables combining the prior
knowledge and the information obtained from measurements. Figs.
\ref{fig_lambdapostfields} and \ref{fig_respposttime} present the
results of the Bayesian inference verification. It is shown that
even for highly nonlinear model, the updating process leads to
much more precise prediction of the material properties as well as
model responses.

The drawback of the described procedure is the high computational
cost. One simulation of the presented experiment takes $3.8$
$\mathrm{[s]}$ at Intel Core Duo Processor T$9600$ with $4$GB RAM.
Hence, the whole sampling lasted almost $90$ $\mathrm{[h]}$.
Therefore, our future work will be focussed on the acceleration of
sampling procedure via approximation of the structural response (see
e.g. \cite{Kucerova:2010:TM}) or approximation of the posterior density
\cite{Rosic:2011:JCP} by polynomial chaos expansion.

\section*{Acknowledgement}
This outcome has been achieved with the financial support of the
Czech Science Foundation, project No. 103/08/1531 and No.
105/11/0411 and the Ministry of Education, Youth and Sports,
project No. MSM6840770003. Many thanks belong also to Prof. H.~G.
Matthies, Ph.D. from TU Braunschweig for a lot of fruitful
discussions through the whole work on this paper.

\bibliographystyle{elsarticle-num}
\bibliography{liter-paper_arxiv}

\begin{thebibliography}{10}
\expandafter\ifx\csname url\endcsname\relax
  \def\url#1{\texttt{#1}}\fi
\expandafter\ifx\csname urlprefix\endcsname\relax\def\urlprefix{URL }\fi
\expandafter\ifx\csname href\endcsname\relax
  \def\href#1#2{#2} \def\path#1{#1}\fi

\bibitem{Tarantola:2005}
A.~Tarantola, Inverse Problem Theory and Methods for Model Parameter
  Estimation, Society for Industrial and Applied Mathematics, 2005.

\bibitem{Kunzel:1997}
H.~K\"{u}nzel, K.~Kiessl, Calculation of heat and moisture transfer in exposed
  building components, International Journal of Heat Mass Transfer 40 (1997)
  159--167.

\bibitem{Ditlevsen:1996}
O.~Ditlevsen, H.~O. Madsen, Structural Reliability Methods, John Wiley \& Sons
  Ltd, Chichester, England, 1996.

\bibitem{Stefanou:2009:CMAME}
G.~Stefanou, The stochastic finite element method: Past, present and future,
  Computer Methods in Applied Mechanics and Engineering 198~(9--12) (2009)
  1031--1051.

\bibitem{Sejnoha:MS:2008}
J.~S\'{y}kora, J.~Vorel, T.~Krej\v{c}\'{\i}, M.~\v{S}ejnoha, J.~\v{S}ejnoha,
  Analysis of coupled heat and moisture transfer in masonry structures,
  Materials and Structures 42~(8) (2009) 1153--1167.

\bibitem{Valenta:2009:IJMCE}
R.~Valenta, M.~\v{S}ejnoha, J.~Zeman, Macroscopic constitutive law for mastic
  asphalt mixtures from multiscale modeling, International Journal for
  Multiscale Computational Engineering 8~(1) (2010) 131--149.

\bibitem{Vorel:2008:SEM}
J.~Vorel, M.~\v{S}ejnoha, Evaluation of homogenized thermal conductivities of
  imperfect carbon-carbon textile composites using the mori-tanaka method,
  Structural Engineering and Mechanics 33~(4) (2009) 429--446.

\bibitem{Ghanem:2003}
R.~Ghanem, P.~D. Spanos, Stochastic finite elements: A spectral approach,
  second revised Edition, Dover Publications, Mineola, New York, 2003.

\bibitem{Kucerova:2007:PHD}
A.~Ku\v{c}erov\'a, Identification of nonlinear mechanical model parameters
  based on softcomputing methods, Ph.D. thesis, Ecole Normale Sup\'erieure de
  Cachan, Laboratoire de M\'ecanique et Technologie (2007).

\bibitem{Lubineau:2009:CM}
G.~Lubineau, A goal-oriented field measurement filtering technique for the
  identification of material model parameters, Computational Mechanics 44~(5)
  (2009) 591--603.

\bibitem{Kucerova:2009:WCSMO}
A.~Ku\v{c}erov\'a, Methodology of field measurements filtering to maximize the
  correlation with material parameters, in: WCSMO-8, Laborat\'orio Nacional de
  Engenharia Civil, Lisboa, Portugal, 2009, pp. on CD--ROM.

\bibitem{Kuraz:2010:JCAM}
M.~Kur\'a\v{z}, P.~Mayer, M.~Lep\v{s}, D.~Trpko\v{s}ov\'a, An adaptive time
  discretization of the classical and the dual porosity model of {R}ichards'
  equation, Journal of Computational and Applied Mathematics 233~(12) (2010)
  3167--3177.

\bibitem{Kucerova:2009:EC}
A.~Ku\v{c}erov\'{a}, D.~Brancherie, A.~Ibrahimbegovi\'{c}, J.~Zeman,
  Z.~Bittnar, Novel anisotropic continuum-discrete damage model capable of
  representing localized failure of massive structures. {P}art {II}:
  identification from tests under heterogeneous stress field, Engineering
  Computations 26~(1/2) (2009) 128--144.

\bibitem{Lehky:2005}
D.~Lehk\'y, D.~Nov\'ak, Probabilistic inverse analysis: Random material
  parameters of reinforced concrete frame, in: Ninth International Conference
  on Engineering Applications of Neural Networks, EAAN2005, Lille, France,
  2005, pp. 147--154.

\bibitem{Kennedy:2001:JRSS}
M.~C. Kennedy, A.~O'Hagan, Bayesian calibration of computer models, Journal of
  the Royal Statistical Society: Series B (Statistical Methodology) 63~(3)
  (2001) 425--464.

\bibitem{Yee:2008:JWEIA}
E.~Yee, F.-S. Lien, A.~Keats, R.~D'Amours, Bayesian inversion of concentration
  data: Source reconstruction in the adjoint representation of atmospheric
  diffusion, Journal of Wind Engineering and Industrial Aerodynamics
  96~(10--11) (2008) 1805--1816.

\bibitem{Fu:2009:JHydro}
J.~Fu, J.~J. G\'omez-Hern\'andez, Uncertainty assessment and data worth in
  groundwater flow and mass transport modeling using a blocking {M}arkov chain
  {M}onte {C}arlo method, Journal of Hydrology 364~(3--4) (2009) 328--341.

\bibitem{Parthasarathy:2008:IJHMT}
S.~Parthasarathy, C.~Balaji, Estimation of parameters in multi-mode heat
  transfer problems using {B}ayesian inference – effect of noise and a priori,
  International Journal of Heat and Mass Transfer 51~(9--10) (2008) 2313--2334.

\bibitem{Tierney:1994:AS}
L.~Tierney, Markov chains for exploring posterior distributions, Annals of
  Statistics 22~(4) (1994) 1701--1728.

\bibitem{Matthies:2007:IB}
H.~G. Matthies, Encyclopedia of Computational Mechanics, John Wiley \& Sons,
  Ltd., 2007, Ch. Uncertainty Quantification with Stochastic Finite Elements.

\bibitem{Rosic:2008:JSSCM}
B.~Rosi\'c, H.~G. Matthies, Computational approaches to inelastic media with
  uncertain parameters, Journal of the Serbian Society for Computational
  Mechanics 2~(1) (2008) 28--43.

\bibitem{Marzouk:2009:JCP}
Y.~Marzouk, H.~Najm, Dimensionality reduction and polynomial chaos acceleration
  of {B}ayesian inference in inverse problems, Journal of Computational Physics
  228~(6) (2009) 1862--1902.

\bibitem{Matthies:2005:CMAME}
H.~Matthies, A.~Keese, Galerkin methods for linear and nonlinear elliptic
  stochastic partial differential equations, Computer Methods in Applied
  Mechanics and Engineering 194~(12-16) (2005) 1295--1331.

\bibitem{Khoromskij:2008}
B.~N. Khoromskij, A.~Litvinenko, Data sparse computation of the
  {K}arhunen-{L}o\`eve expansion, in: Proceedings of International Conference
  on Numerical Analysis and Applied Mathematics 2008, Vol. 1048, 2008, pp.
  311--314.

\bibitem{Cerny:2009:CMEM}
R.~\v{C}ern\'y, J.~Mad\v{e}ra, J.~Ko\v{c}\'i, E.~Vejmelkov\'a, Heat and
  moisture transport in porous materials involving cyclic wetting and drying,
  in: Computational Methods and Experimental Measurements XIV, Vol.~48 of WIT
  Transactions on Modelling and Simulation, 2009, pp. 3--12.

\bibitem{Pavlik:SFR:2010}
Z.~Pavl\'{i}k, J.~Mihulka, J.~\v{Z}um\'{a}r, R.~\v{C}ern\'{y}, Experimental
  monitoring of moisture transfer across interfaces in brick masonry, in:
  Structural Faults and Repair, 2010.

\bibitem{Lombardo:2009:IJMCE}
M.~Lombardo, J.~Zeman, M.~\v{S}ejnoha, G.~Falsone, Stochastic modeling of
  chaotic masonry via mesostructural characterization, International Journal
  for Multiscale Computational Engineering 7~(2) (2009) 171--185.

\bibitem{Kucerova:2010:TM}
A.~Ku\v{c}erov\'a, H.~G. Matthies, Uncertainty updating in the description of
  heterogeneous materials, Technische Mechanik 30~(1--3) (2010) 211--226.

\bibitem{Rosic:2011:JCP}
B.~Rosi\'{c}, A.~Litvinenko, O.~Pajonk, H.~G. Matthies, Direct {B}ayesian
  update of polynomial chaos representations, Journal of Computational Physics
  0~(0) (2011) Submitted for publication.

\end{thebibliography}
\end{document}